\documentclass{amsart}

\usepackage{amsmath}
\usepackage{amssymb}
\usepackage{amsfonts}
\usepackage{MnSymbol}
\usepackage{graphics} 
\usepackage{epsfig} 
\usepackage{psfrag}
\usepackage{pdfsync}
\usepackage[usenames]{color}
\usepackage{cancel}
\usepackage[pagewise]{lineno}

\definecolor{arancio}{rgb}{0.90,0.50,0.20}

\definecolor{blu}{rgb}{0.,0.,1.}

\definecolor{pavone}{rgb}{0.00,0.00,0.63}

\definecolor{malva}{rgb}{0.10,0.50,0.50}

\definecolor{rosso}{rgb}{1.,0.,0.}

\definecolor{geranio}{rgb}{0.90,0.00,0.20}

\definecolor{cerulean}{rgb}
{0.0, 0.48, 0.65}

\newtheorem{theorem}{Theorem}[section]
\newtheorem{corollary}[theorem]{Corollary}

\newtheorem{lem}[theorem]{Lemma}
\newtheorem{prop}[theorem]{Proposition}

\theoremstyle{definition}
\newtheorem{definition}{Definition}[section]
\newtheorem{remark}{Remark}[section]

\newcommand{\ep}{\epsilon}
\newcommand{\N}{\mathbb{N}}
\newcommand{\R}{\mathbb{R}}
\newcommand{\lla}{\left\langle}
\newcommand{\rra}{\right\rangle}
\date{\today}

\newcommand{\bcl}{\begin{center}}
\newcommand{\ecl}{\end{center}}
\newcommand{\brl}{\begin{right}}
\newcommand{\erl}{\end{right}}
\newcommand{\ben}{\begin{enumerate}}
\newcommand{\barr}{\begin{array}}
\newcommand{\earr}{\end{array}}
\newcommand{\btab}{\begin{tabular}}
\newcommand{\etab}{\end{tabular}}
\newcommand{\bdoc}{\begin{document}}
\newcommand{\edoc}{\end{document}}
\newcommand{\beqy}{\begin{eqnarray}}

\newcommand{\beq}{\begin{equation}}
\newcommand{\beqi}{\begin{eqnarray*}}
\newcommand{\bitem}{\begin{itemize}}
\newcommand{\brem}{\begin{remark}}
\newcommand{\erem}{\end{remark}}
\newcommand{\eitem}{\end{itemize}}
\newcommand{\nln}{\newline}
\newcommand{\newt}{\newtheorem}

\renewcommand{\a }{\alpha }
\renewcommand{\b }{\beta }
\newcommand{\g }{\gamma}
\newcommand{\G }{\Gamma }
\renewcommand{\d }{\delta }

\newcommand{\D }{\Delta }
\newcommand{\e }{\epsilon }
\newcommand{\z }{\zeta }
\renewcommand{\l }{\lambda }
\renewcommand{\L }{\Lambda }
\newcommand{\m }{\mu }
\newcommand{\n }{\tau }
\renewcommand{\r }{\rho }
\newcommand{\s }{\sigma }
\newcommand{\Sig }{\Sigma }
\renewcommand{\t }{\tau }
\newcommand{\var }{H }
\renewcommand{\o }{\omega }
\renewcommand{\O }{\Omega }

\newcommand{\supp}{\text{\rm supp}\,}
\newcommand{\sgn}{\text{\rm sgn}\,}

\title[Signed Radon measure-valued solutions]
{Signed Radon measure-valued solutions \\ of flux saturated scalar conservation laws}

\author[Bertsch]{Michiel Bertsch}
\address{Dipartimento di Matematica, Universit\`a di Roma ``Tor Vergata'', 
Via della Ricerca Scientifica, 00133 Roma, Italy \\ and
Istituto per le Applicazioni del Calcolo "M. Picone", CNR, Roma, Italy} 
\email{bertsch.michiel@gmail.com}
\author[Smarrazzo]{Flavia Smarrazzo}
\address{Universit\`a Campus Bio-Medico di Roma\\ Via Alvaro del Portillo 21, 00128 Roma, Italy}
\thanks{}
\email{flavia.smarrazzo@gmail.com}
\author[Terracina]{Andrea Terracina}
\address{Dipartimento di Matematica "G. Castelnuovo", Universit\`a ``Sapienza'' di Roma\\ P.le A. Moro 5, I-00185 Roma, Italy}
\email{terracina@mat.uniroma1.it}
\author[Tesei]{Alberto Tesei}
\address{Dipartimento di Matematica "G. Castelnuovo", Universit\`a ``Sapienza'' di Roma\\ P.le A. Moro 5, I-00185 Roma, Italy, and 
Istituto per le Applicazioni del Calcolo "M. Picone", CNR, Roma, Italy}
\email{albertotesei@gmail.com}

\subjclass{Primary:  Secondary: 
 }  


\date{\today}

\begin{document}

\bibliographystyle{h-elsevier2}
\begin{abstract} 

We prove existence and uniqueness  for a class of signed Radon measure-valued entropy solutions of the Cauchy problem for a first order scalar hyperbolic conservation law in one space dimension. The initial data of the problem is a finite superposition of Dirac masses, whereas the flux is Lipschitz continuous and bounded. The solution class is determined by an additional condition which is needed to prove uniqueness.
\end{abstract}

\maketitle


\section{Introduction}\label{intro}

We study the Cauchy problem for the scalar conservation law:
$$
\left\{\begin{array}{ll}
u_t+\left[H(u)\right]_x=0 
& \quad\mbox{in}\   \R\times (0,T)=:S \smallskip\\
u=u_0  &\quad\mbox{in}\  \R\times \{0\}\,,
\end{array}\right. \leqno{(CL)}
$$ 
where 
$$
H\in  W^{1,\infty}(\R)\,,\quad H(0)=0 \leqno{(A_0)}
$$
(obviously the condition $H(0)=0$ is not restrictive).
The initial condition $u_0$ is a signed Radon measure on $\R$.
In most of the paper we shall assume that its singular part, $u_{0s}$,  
is a finite superposition of Dirac masses:
$$
u_{0s}=\sum_{j=1}^{p} c_j \delta_{x_j}\qquad 
\text{($x_1<x_2<\dots<x_p$; \ $c_j\in \R\setminus\{0\}$ for $1\le j\le p$)\,.}
\leqno{(A_1)}
$$
In that case we denote  the support of the singular measure $u_{0s}$ by $F$:
$$
F=\{x_1,x_2,\cdots,x_p\}.
$$

In \cite{BSTT2} we considered the case of {\it nonnegative} initial measures $u_0$.
In the present paper we consider the case of {\it signed} measures (see \cite{BSTT1, DS, F, LP} 
for motivations and related remarks). A specific motivation is the link between 
measure-valued solutions of $(CL)$ and
{\em discontinuous} solutions of the Cauchy problem for the Hamilton-Jacobi equation
$$
\left\{\begin{array}{ll}
U_t+H(U_x)=0 
& \quad\mbox{in}\ S \smallskip\\
U=U_0  &\quad\mbox{in}\  \R\times \{0\}\,,
\end{array}\right. \leqno{(HJ)}
$$ 
where $U_0\in BV_{loc}(\R)\cap L^\infty(\R)$, 
$U_0'\in L^1_{loc}(\R\setminus F)$, 
and
$U_0(x_j^+)\neq U_0(x_j^-)$ if  $x_j\in F$. 
If $(A_1)$ is satisfied, the distributional derivative $U_0'$ is a Radon measure without singular continuous part,
$U_0'=\sum\limits_{j=1}^p\left[U_0(x_j^+)-U_0(x_j^-)\right] \delta_{x_j} + (U_0')_{ac}$,
and problems $(CL)$, $(HJ)$ are formally related by the equality $u=U_x$. In a forthcoming paper \cite{BSTT4}, problem $(HJ)$ will  
be studied in the context of viscosity solutions.

It is known (\cite{BSTT1}) that $(i)$ the singular part $u_s$ of a suitably defined entropy solution 
may persist for some positive time (see \cite[Theorem 3.5]{BSTT1}) and $(ii)$ 
entropy solutions are not always 
uniquely determined by the initial condition $u_0$ (see also Remark \ref{nonuniq}). 
To overcome the latter problem, we introduced in \cite{BSTT2}
a so-called compatibility condition at those points where $u_s(\cdot,t)$ is a Dirac mass, 
and used it 
as a uniqueness criterion for nonnegative measure-valued solutions.  

The starting point of the present paper is the statement that, for general signed initial 
measures $u_0$, the singular part $u_s$ of any local entropy solution $u$ of $(CL)$ (in the sense of  
Definition~\ref{enso}) satisfies a monotonicity result: both the 
positive and negative part of $u_s$, $[u(\cdot,t)]_s^\pm$, 
are nonincreasing with respect to $t$ (see Theorem~\ref{recl}). 
For the class of initial measures satisfying $(A_1)$ this implies that 
the support of the singular part $u_s$ of an entropy solution of problem $(CL)$ 
is a subset of $F\times [0,T]$ and, in addition, that the sign of $[u(\cdot,t)]_s$ is determined 
by that of $u_{0s}$. Having this in mind it is rather straightforward to adapt the concept of 
compatibility condition in \cite{BSTT2} to {\it signed} measure-valued solutions (see Definition 
\ref{compatibility condition}).

The main result of the paper is that  if $(A_0)$ and $(A_1)$ are satisfied, then  $(CL)$ is well-posed
in the class of entropy solutions which satisfy the compatibility condition at the $p$ points $x_j\in F$.

Existence of a solution is proven by a constructive approach which can be outlined as follows.
By $(A_0)$-$(A_1)$ 
there exists a positive time $\tau$ until which all singularities persist (see \cite[Theorem 3.5]{BSTT1}), thus the real line is the disjoint union of $p+1$ intervals. In each interval we solve the initial-boundary value problem for the conservation law in $(CL)$, the initial data being the restriction of $u_{0r}$ to that interval, with ``boundary conditions equal to infinity". Namely, we consider the {\em singular} 
Dirichlet  initial-boundary value problems 
\begin{equation}\label{dj}
\left\{\begin{array}{ll}
u_{t}+[H(u)]_x=0&\mbox{in}\ \,(x_{j-1},x_j)\times(0,T) \\
u=\pm\infty&\mbox{in $\{x_{j-1}\}\times (0,T)$} \\
u=\pm\infty&\mbox{in $\{x_j\}\times (0,T)$} \\
u=u_{0r}&\mbox{in $(x_{j-1},x_j)\times\{0\}$}\,
\end{array}\right.
\end{equation}
with $j=2,\dots,p$,  and
\begin{equation}\label{d1}
\left\{\begin{array}{ll}
u_{t}+[H(u)]_x=0&\mbox{in}\ \,(-\infty,x_1)\times(0,T) \\
u=\pm\infty&\mbox{in $\{x_1\}\times (0,T)$} \\
u=u_{0r}&\mbox{in $(-\infty,x_1)\times\{0\}$}\,,
\end{array}\right.
\end{equation}
\begin{equation}\label{dp}
\left\{\begin{array}{ll}
u_{t}+[H(u)]_x=0&\mbox{in}\ \,(x_p,\infty)\times(0,T) \\
u=\pm\infty&\mbox{in $\{x_p\}\times (0,T)$} \\
u=u_{0r}&\mbox{in $(x_p,\infty)\times\{0\}$}\,.
\end{array}\right.
\end{equation}
The choice between $u=\infty$ and $u=-\infty$ at $x_j$ is determined by the sign of $c_j$: we choose
$\infty$ if $c_j>0$ and $-\infty$ if $c_j<0$.
Existence and uniqueness of an entropy solution to each problem \eqref{dj}-\eqref{dp} is proven in Sections \ref{unid}-\ref{exid}. In particular, existence follows from an approximation procedure which makes use of BV initial and boundary data, avoiding the $L^\infty$-theory of initial-boundary value problems developed in \cite{O} (see Section \ref{exid}).

The function determined  by solutions of  \eqref{dj}-\eqref{dp} in $\R\times(0,\t)$ is, by definition, the regular part of a Radon measure, whose singular part is defined by observing that the variation of mass at each point $x_j$ depends on the sweeping effect of the flux across $x_j$ (see \eqref{C_j+}, \eqref{C_j-}, \eqref{defus} and Proposition \ref{letra}).
Then it is proven that this measure is the unique entropy solution of $(CL)$ (in the sense of Definition \ref{enso}) which satisfies the compatibility conditions at all $x_j \in F$ until the time $t=\tau$. 
Here we use that the required compatibility condition 
for the solution of the Cauchy problem $(CL)$ at $x_j$  
is exactly the entropic formulation of the boundary conditions "$u=\pm\infty$" for the 
singular Dirichlet problems (see also Remark~\ref{comp-entr}).
If $\tau<T$ we iterate the procedure in $\R\times(\tau,T)$ with a smaller number a singularities, thus well-posedness of $(CL)$ follows in a finite number of steps (see Section \ref{cleu}). We observe that the proof of uniqueness of entropy solutions to problem $(CL)$ relies on a general comparison principle between entropy sub and super-solutions of \eqref{dj}--\eqref{dp} (see Definitions \ref{defsub}--\ref{-infb} and Theorem \ref{cmp} below) which is independent of the above construction procedure. In this sense the comparison results are stronger than those in \cite[Theorem 3.2]{BSTT2}.

The results in the paper can be esaily extended to the case that $u_{0s}$ is a locally finite superposition of Dirac masses (namely, if the number of Dirac masses in every bounded interval is finite).


\section{Preliminaries}\label{preli}
\setcounter{equation}{0}

Let $\chi_E$ denote the characteristic function of $E\subseteq\R$. For every $u\in\R$ we set
\begin{equation*}
[u]_\pm:=\max\{\pm u,0\}, \quad{\rm sgn}_\pm(u):=\pm\chi_{\R_\pm}(u), \quad {\rm sgn}(u):={\rm sgn}_-(u)+{\rm sgn}_+(u)\,.
\end{equation*}
For every real function $f$ on $\R$ and $x_0\in\R$ we say that 
$$
\text{{\rm ess\,}$\lim_{x\to x_0^{\pm}}  f(x) =l\in\R$\,,}
$$ 
if there is a  null set $E^*\subseteq\R$  such that $f(x_n)\to l$ if $\left\{x_n\right\}\subseteq\R\setminus\! (E^*\!\cup\!\{x_0\})$, $x_n\to x_0^{\pm}$. 

For every open subset $\Omega\subseteq\R$ we denote by $C_c(\Omega)$ the space  of continuous real functions with compact support in $\Omega$ 
and by $\mathcal{M}^+(\Omega)$ the cone of the nonnegative Radon measures on $\Omega$. According to  \cite[Section 1.3]{EG}, we say that $\nu$ is a (signed) Radon measure on $\Omega$ if there exist a (nonnegative) Radon measure $\mu\in\mathcal{M}^+(\Omega)$ and a locally $\mu$-summable function $f:\Omega \to[-\infty,\infty ]$ such that $$\nu(K)=\int_Kf\,d\mu$$ for all compact sets $K\subset \Omega$. The space of (signed) Radon measures on $\Omega$ will be denoted by $\mathcal{M}(\Omega)$.   

If $\mu,\nu\in\mathcal{M}(\Omega)$, we say that $\mu\le\nu$ in $\mathcal{M}(\Omega)$ if $\nu-\mu\in\mathcal{M}^+(\Omega)$. We denote by $\lla \cdot, \cdot\rra_{\Omega}$ the duality map between $\mathcal M(\Omega)$  and  $C_c(\Omega)$. A sequence $\{\mu_n\}$ of Radon measures on $\R$ converges weakly* to a Radon measure $\mu$, $\mu_n\stackrel{*}\rightharpoonup \mu$, if $\left\langle \mu_n,\rho\right\rangle_{\R}\to \left\langle \mu,\rho\right\rangle_{\R}$ for all $\rho \in C_c(\R)$. For any compact $K\in\R$ the space $\mathcal{M}(K)$ is a Banach space with norm $\|\mu \|_{\mathcal M(K)}:=|\mu|(K)$,
where $|\m|$ denotes the total variation of $\m$.
A sequence $\{\mu_n\}$ converges strongly to $\mu$ in $\mathcal{M}(K)$ if  $\|\mu_n-\mu\|_{\mathcal{M}(K)}\to 0$ as  $n\to\infty$.
Similar definitions are used for Radon measures on any subset of $S:=\R\times (0,T)$. 

Every $\mu\in\mathcal{M}(\R)$ has a 
unique decomposition $\mu=\mu_{ac}+\mu_s$, with $\mu_{ac} \in\mathcal{M}(\R)$
absolutely continuous and $\mu_s\in\mathcal{M}(\R)$ singular with respect to the Lebesgue  measure. 
We  denote by $\mu_r\in L^1_{loc}(\R)$ the density of $\mu_{ac}$. Every function $f\in L^1_{loc}(\R)$ can be identified to an absolutely continuous Radon measure on $\R$; we shall denote this measure by the same symbol $f$ used for the function. 

The restriction $\mu \, \lefthalfcup E$ of $\mu\in \mathcal{M}(\R)$ to a Borel set $E\subseteq \R$ is defined by 
$(\mu \, \lefthalfcup E) (A):=\mu(E\cap A)$ for any Borel set $A\subseteq\R$. 
Similar notations are used for $\mathcal{M}(S)$.

For every open subset $\Omega\subseteq\R$ we denote by $BV(\Omega)$ the Banach space of functions of bounded variation in $\Omega$: 
\begin{equation*}
BV(\Omega):= \left\{z\in L^1(\Omega) \, | \, z' \in \mathcal{M}(\Omega), \|z'\|_{\mathcal{M}(\Omega)}<\infty \right\}\,,\quad
\|z\|_{BV(\Omega)}:= \|z\|_{L^1(\Omega)}+ \|z'\|_{\mathcal{M}(\Omega)}\,,
\end{equation*}
where $z'$ is the first order distributional derivative. The total variation in $\O$ of $z$ is $TV(z;\O):=\|z'\|_{\mathcal{M}(\Omega)}$. We say that $z\in BV_{loc}(\R)$ if $z\in BV(\Omega)$ for every open subset $\Omega\subset\subset\R$.

In the remainder of this section $\Omega$ denotes an open subset of $\R$, and 
$Q_T=\Omega \times (0,T)$. By $C([0,T];\mathcal{M}(\Omega))$ we denote the subset of strongly continuous mappings from $[0,T]$ into $\mathcal{M}(\Omega)$ - namely, $u\in C([0,T];\mathcal{M}(\Omega))$ if for all $t_0\in[0,T]$ and for every compact $K\in\Omega$ there holds $ \|u(\cdot,t)-u(\cdot,t_0)\|_{\mathcal{M}(K)}\to0$ as $t\to t_0$.

\begin{definition}\label{dli}
We denote by $L^{\infty}(0,T;\mathcal{M}^+(\Omega))$ the set of nonnegative Radon measures $u\in \mathcal{M}^+(Q_T)$ such that 
for a.e.~$t\in (0,T)$ there is a measure $u(\cdot,t)\in \mathcal{M}^+(\Omega)$ with the following properties:

\noindent $(i)$ if $\zeta\in C([0,T];C_c(\Omega))$ the map $t\mapsto \left\langle u(\cdot,t),\zeta(\cdot,t)\right\rangle_{\Omega}$ 
belongs to $L^1(0,T)$ and
\begin{equation}\label{eq.disintegrazioneU}
\left\langle u,\zeta\right\rangle_{Q_T}=\int_0^T\left\langle u(\cdot,t),\zeta(\cdot,t)\right\rangle_{\Omega}\,dt\,;
\end{equation}

\noindent $(ii)$ the map $t\mapsto \|u(\cdot,t)\|_{\mathcal{M}(K)}$ belongs to $L^{\infty}(0,T)$ for every compact $K\subset\Omega$.
\end{definition} 

\begin{remark} Definition \ref{dli} implies that for all  $\rho\in C_c(\Omega)$ the map  $t\mapsto \left\langle u(\cdot,t),\rho\right\rangle_{\Omega}$ 
is measurable, thus the map   $u:(0,T) \to \mathcal{M}^+(\Omega)$ is weakly* measurable. For simplicity we prefer the notation $L^{\infty}(0,T;\mathcal{M}^+(\Omega))$ to the more correct one $L^{\infty}_{w*}(0,T;
\mathcal{M}^+(\Omega))$. Moreover, as a consequence of Definition \ref{dli}-$(i)$, it can be seen that for every Borel set $E\subseteq Q_T$ the map $t\to u(\cdot,t)\big(E^t\big)$ is Lebesgue measurable and there holds
\begin{equation}
u(E)=\int_0^Tu(\cdot,t)\big(E^t\big)\,dt\,\qquad(E^t=\{x\in\Omega\,:\ \,(x,t)\in E\}).
\end{equation}\end{remark}

If $u\in L^{\infty}(0,T;\mathcal{M}^+(\Omega))$, then $u_{ac},\, u_s\in L^{\infty}(0,T;\mathcal{M}^+(\Omega))$ as well, and  $u_r\in L^{\infty}(0,T;L^1_{\rm loc}(\Omega))$.
 Moreover, equality \eqref{eq.disintegrazioneU} implies 
$$\mbox{$\lla u_{ac}, \zeta\rra_{Q_T}=\int\!\!\!\int_{Q_T}u_r \,\zeta\,dxdt\ \ $ and $\ \ 
\lla u_s, \zeta \rra_{Q_T} =\int_0^T\!\! \lla u_s(\cdot,t),\zeta(\cdot,t) \rra_{\Omega} dt$.}$$ 
Denoting by $[u(\cdot,t)]_{ac},\,[u(\cdot,t)]_s\in\mathcal{M}^+(\Omega)$ the absolutely continuous and the singular part of the measure $u(\cdot,t)\in \mathcal{M}^+(\Omega)$, a routine proof shows that for $a.e.$ $t\in (0,T)$ 
\begin{equation}\label{us(t)=u(t)s}
u_s(\cdot,t)=[u(\cdot,t)]_s\,,\quad u_{ac}(\cdot,t)=[u(\cdot,t)]_{ac}\,, \quad u_r(\cdot,t)=[u(\cdot,t)]_r\,,
\end{equation}
where $[u(\cdot,t)]_r$ denotes the density of the measure $[u(\cdot,t)]_{ac}$. In view of \eqref{us(t)=u(t)s}, we shall always  identify the quantities which appear 
on either side of equalities \eqref{us(t)=u(t)s}. 
\smallskip

We say that a (signed) Radon measure $u\in\mathcal{M}(Q_T)$ belongs to $L^{\infty}(0,T;\mathcal{M}(\Omega))$ if both  $u^+$ and $u^-$ belong to $L^{\infty}(0,T;\mathcal{M}^+(\Omega))$. In particular, this implies that:
\smallskip

\noindent $(\alpha)$ the total variation $|u|$ of the measure $u$ belongs to $L^{\infty}(0,T;\mathcal{M}^+(\Omega))$;

\noindent $(\beta)$ conditions $(i)$ and $(ii)$ of Definition \ref{dli} hold with $u(\cdot,t):=u^+(\cdot,t)-u^-(\cdot,t)$ for $a.e.\ t\in (0,T)$.  
\smallskip

Moreover, since $u^+$ and $u^-$ are mutually singular, it follows that for $a.e.$ $t$ the nonnegative measures $u^+(\cdot,t)$ and $u^-(\cdot,t)$ are mutually singular, whence
\begin{equation}\label{eq.scambio.pm}
u^{\pm}(\cdot,t)=[u(\cdot,t)]^{\pm}\,,\quad \ |u(\cdot,t)|=|u|(\cdot,t)\quad\mbox{for}\ \,a.e.\ \,t\in (0,T)\,,
\end{equation} 
and
\begin{equation}\label{uspm}
u_s^{\pm}(\cdot,t)=[u(\cdot,t)]_{s}^{\pm},\quad |u_s|(\cdot,t)=\left|[u(\cdot,t)]_{s}\right|\quad\mbox{for}\ \,a.e.\ \,t\in (0,T)\,.
\end{equation}


\section{Results}\label{resu}
\setcounter{equation}{0}

\noindent For any $\tau\in(0,T]$ and open subset $\Omega\subseteq\R$ set $Q_{\tau}:=\Omega\times (0,\tau]$, $Q_T\equiv Q$; set also $S_{\tau}:=\R\times (0,\tau]$, $S_T\equiv S$. Solutions of problem $(CL)$ are meant in the following sense. 

\begin{definition}\label{deso}  
Let $u_0$ be a signed Radon measure on $\Omega$ and let $(A_0)$ be satisfied.
A measure $u\in L^{\infty}(0,T;\mathcal{M}(\Omega))$ is a {\it solution} of problem $(CL)$ in $Q_{\tau}$  
if for all $\zeta\in C^1([0,\tau];C^1_c(\Omega))$, $\zeta(\cdot,\tau)=0$ in $\Omega$ there holds
\begin{equation}\label{ewf}
\iint_{Q_{\tau}} \big[u_r\zeta_t+ H(u_r) \,\zeta_x\big]\,dxdt+\int_0^{\tau}\lla u_s(\cdot,t), 
\zeta_t(\cdot,t)\rra_{\Omega}dt= - \lla u_0,\zeta(\cdot,0)\rra_{\Omega}\,.
\end{equation}
Solutions of $(CL)$ in $S$ are simply referred to as ``solutions of $(CL)$''.
\end{definition}

\begin{definition}\label{enso} Let $u_0$ be a signed Radon measure on $\Omega$ and let $(A_0)$ be satisfied. A solution  of $(CL)$ in $Q_{\tau}$ is called an {\em entropy solution in $Q_{\tau}$} if 
it satisfies the {\em entropy inequality}
\begin{eqnarray}\label{mkru}
&&\iint_{Q_{\tau}} \left\{|u_r-k|\,\zeta_t+\sgn(u_r-k)\left [H(u_r)-H(k)\right ]\zeta_x\right\}dxdt 
\,+\\
&&\quad + \int_0^\tau\left\langle |u_s(\cdot,t)|,\zeta_t(\cdot,t)\right\rangle_{\Omega}\,dt \geq - 
\int_{\Omega} |u_{0r}(x)-k|\,\zeta(x,0)
\,dx - \left\langle |u_{0s}|, \zeta(\cdot,0)\right\rangle_{\Omega} \nonumber  
\end{eqnarray}
 for all $\zeta\in C^1([0,\tau];C^1_c(\Omega))$, $\zeta\geq0$, $\zeta(\cdot,\tau)=0$ in $\Omega$, and for all $k\in\R$;
\end{definition}

If $Q_\tau\ne Q_T$, an (entropy) solution in $Q_\tau$ can be considered as a {\it local} (entropy)
solution of $(CL)$. 
For general initial measures, local entropy solutions satisfy the following monotonicity result. 
 \begin{theorem}\label{recl}
Let $(A_0)$ be satisfied, let $u_0$ be a signed Radon measure on $\Omega$ and let $u$
be an entropy solution $u$ of problem $(CL)$ in $Q_T$. 
Then, for a.e.~$0< t_1<t_2< T$, 
there holds
\begin{equation}\label{mono.us}
 [u(\cdot,t_2) ]_s^{\pm} \leq [u(\cdot,t_1) ]_s^{\pm}\le u_{0s}^{\pm} \quad\mbox{in}\ \,\mathcal{M}(\Omega)\,.
\end{equation}
\end{theorem}

Now we consider the case that $u_{0s}$ is the sum of a finite number of Dirac masses with 
support $F$.

\begin{remark}\label{continuity}
Let $(A_0)-(A_1)$ be satisfied and let $u$ be an entropy solution of problem $(CL)$ in $Q_T$.
Arguing as in the proof of Proposition 3.20 in \cite{BSTT1}, it follows that 
$u\in C((0,T];\mathcal{M}(\Omega))$.
\end{remark}

\begin{corollary}\label{recl-bis}
Let $(A_0)-(A_1)$ be satisfied and let $u$
be an entropy solution $u$ of problem $(CL)$ in $Q_T$.
Then $u_s\in C([0,T];\mathcal{M}(\Omega))$, $[u(\cdot,0)]_s=u_{0s}$, 
\eqref{mono.us} holds for any $0\le t_1\le t_2 \le T$ and
for every $x_j\in F\cap \Omega$ there exists  $t_j\in(0,T]$ such that
  \begin{equation}\label{dewait}
 \begin{cases}
u_s(\cdot,t)(\{x_j\})\neq 0 &\text{if $t\in[0,t_j)$,} \\
u_s(\cdot,t)(\{x_j\})=0 &\text{if $t\in(t_j,T ]$.} 
\end{cases} 
\end{equation}
\end{corollary}

\noindent We observe that the proof of Corollary \ref{recl-bis} provides an explicit lower bound for $t_j$. 
\smallskip

If $x_j\in F$ and $t_j\in (0,T]$ as in Corollary \ref{recl-bis}, Theorem \ref{recl} 
implies that the support of the singular part of any entropy solution is a subset of  $F\times [0,T]$ and that
the Delta mass at $x_j\in F$ does not change sign in the interval $[0,t_j)$.
Therefore we may formulate a compatibility condition at $x_j$ which 
depends on the sign of $c_j$, i.e.~on the sign of the initial Delta mass at $x_j$:

\begin{definition}\label{compatibility condition} Let $(A_0)-(A_1)$ be satisfied. An entropy solution $u$ of $(CL)$ in $Q_T$ is said to satisfy the 
{\bf compatibility condition} at $x_j\in F\cap \Omega$ if 
\begin{subequations}
\begin{equation}\label{comp1}
{\rm ess} \lim_{x\to x_j^+} \int_0^{t_j}\!{\sgn}_\pm(u_r(x,t)-k)\, \big[H(u_r(x,t))-H(k)\big]\,\beta(t)\,dt  \le 0
\quad \text{if }\pm c_j<0,
\end{equation}
\begin{equation}\label{comp2}
{\rm ess} \lim_{x\to x_j^-} \int_0^{t_j}\!{\sgn}_\pm(u_r(x,t)-k)\, \big[H(u_r(x,t))-H(k)\big]\,\beta(t)\,dt\ge0 
\quad \text{if }\pm c_j<0
\end{equation}
\end{subequations}
for all $\beta\in C^1_c(0,t_j)$, $\beta\ge 0$ and $k\in\R$, where $t_j\in (0,T]$ is defined by 
Corollary \ref{recl-bis}.
\end{definition}

We shall prove below (see Remark \ref{conneHJ}) that, if $(A_0)$-$(A_1)$ hold, 
for every entropy solution $u$ of $(CL)$ the limits
\begin{equation}
{\rm ess} \lim_{x\to x_j^{\pm}} \int_0^T {\sgn}_\pm(u_r(x,t)-k)[H(u_r(x,t))-H(k)]\,\beta(t)\,dt \qquad ( j=1,\dots, p
)\,, 
\end{equation}
with $\b$ and $k$ as above, exist and are finite. Hence Definition 
\ref{compatibility condition} is well-posed.

\smallskip

The main result of the paper is the well-posedness of problem $(CL)$, if $u_0$ satisfies $(A_1)$, 
in the class of entropy solutions in $C([0,T];\mathcal{M}(\R))$ which satisfy the 
compatibility condition in $\supp u_0$.

\begin{theorem}\label{exiuni}
 Let $(A_0)$-$(A_1)$ be satisfied. Then there exists  a unique entropy solution of problem $(CL)$
  which belongs to $C([0,T];\mathcal{M}(\R))$ and satisfies the compatibility condition at 
 all $x_j\in \supp u_0$.
 \end{theorem}

\begin{remark}\label{nonuniq} It was already observed in \cite{BSTT1} that in general measure-valued
entropy solutions are not unique. This is essentially a consequence of the elementary 
observation that there
exists a unique entropy solution for which $[u_s(t)]=u_{0s}$ for a.e.~$t\in (0,T)$ (it is enough
to set $u=u_0+\tilde u$, where $\tilde u$ is the entropy solution with initial data $u_{0r}$).
But if $u_0$ satisfies $(A_1)$ and $F\ne \emptyset$, one easily checks that  
if the function $H$, satisfying $(A_0)$,  is not constant in intervals of the type $(a,\infty)$ and $(-\infty,b)$, 
then such solution does not satisfy  the compatibility condition at $x_j\in F$. In particular,
it does not coincide with  the solution defined by Theorem~\ref{exiuni}.
\end{remark}


\section{Monotonicity of $u_s$.}

In this section we prove Theorem \ref{recl} and Corollary  \ref{recl-bis}.
\medskip

\noindent {\it Proof of Theorem \ref{recl}.} By \eqref{ewf}, for every $k\in\R$ we get
\begin{eqnarray*}
&&\iint_{S} (u_r-k)\beta'(t)\rho(x)\,dxdt + \int_0^T\left\langle u_s(\cdot,t),\rho\right\rangle_{\R}\beta'(t)\,dt 
+\\&&\qquad 
+ \iint_S \left [H(u_r)-H(k) \right ]\,\rho'(x)\beta(t)\,dxdt
=\nonumber\\&&\qquad 
= - \beta(0)\left\{\int_{\R}(u_{0r}-k)\rho(x)\,dx +\left\langle u_{0s}\,,\,\rho\right\rangle_{\R}\right\}\nonumber
\end{eqnarray*}
for all $\rho\in C^1_c(\Omega)$ and $\beta \in C^1_c([0,T))$. By summing and subtracting the above equality from the entropy inequality \eqref{mkru}, for every nonnegative $\rho$ and $\beta$ as above  we obtain
\begin{eqnarray*}
&&\iint_{S} [u_r-k]_\pm\,\beta'(t)\rho(x)\,dxdt + \int_0^T\left\langle [u(\cdot,t)]_s^\pm,\rho\right\rangle_{\R}\beta'(t)\,dt +\\
&&\qquad + \iint_S {\rm sgn}_\pm(u_r-k)\left [H(u_r)-H(k) \right ]\,\rho'(x)\beta(t)\,dxdt\geq \nonumber\\
&&\qquad \geq- \beta(0)\left\{\int_{\R}[u_{0r}-k]_\pm\,\rho(x)\,dx +\left\langle u_{0s}^\pm\,,\,\rho\right\rangle_{\R}\right\}\,,\nonumber
\end{eqnarray*}
Letting $k\to \infty$ with "+" and $k\to -\infty$ with "-", we obtain that 
$$
\int_0^T\left\langle [u(\cdot,t)]_s^\pm,\rho\right\rangle_{\R}\beta'(t)\,dt
\geq
\left\langle u_{0s}^\pm\,,\,\rho\right\rangle_{\R}.
$$
Let $0<t_1<t_2\leq T$. By standard approximation arguments we can choose  
$$
\beta(t)=\beta_n(t)=n(t-t_1)\chi_{[t_1, t_1+1/n]}(t)+\chi_{(t_1+1/n,t_2-1/n)}(t)+n
(t_2-t)\chi_{[t_2-1/n, t_2]}\,.
$$
Arguing as in the proof of Proposition 3.8$(i)$ in \cite{BSTT1}, there exists a null set $N\in (0,T)$
which does not depend on the function $\rho$ such that, letting $n\to \infty$,
\begin{equation}\label{monot4}
\left\langle [u(\cdot,t_2)]_s^\pm\,,\,\rho\right\rangle_{\R} \leq
\left\langle [u(\cdot,t_1)]_s^\pm\,,\,\rho\right\rangle_{\R} \quad\text{if $t_1,t_2\not \in N$}\,.
\end{equation}
Hence the first inequality in \eqref{mono.us} follows from the arbitrariness of $\rho$. 

The second inequality in \eqref{mono.us} can be proved in a similar way, replacing $\beta_n$ by
$$
\beta_n(t)=\chi_{(0,t_1-1/n)}(t)+n
(t_1-t)\chi_{[t_1-1/n, t_1]}\,.
$$
\hfill$\square$

\noindent {\it Proof of Corollary  \ref{recl-bis}.}  Arguing as in the proof of Theorem \ref{recl}, 
for every $\rho \in C^1_c(\Omega)$  and  $t \in (0, T]$ from \eqref{ewf} we get
\begin{eqnarray}\label{w1}
&&\left\langle u(\cdot,t)\,,\,\rho\right\rangle_{\R}- \left\langle u_0\,,\,\rho\right\rangle_{\R} =  \int^{t}_{0}\!\!\!\int_{\R} H(u_r)\,\rho'(x)\,dxds \,.
\end{eqnarray}
Fix any $x_j\in F\cap\Omega$. By standard approximation arguments we can choose in \eqref{w1}
$$
\rho(x)=\rho_n(x)=n(x-x_j+1/n)\chi_{[x_j-1/n, x_j]}+n(x_j+1/n-x)\chi_{(x_j, x_j+1/n]}\,.
$$
Then letting $n\to \infty$, and observing that
\begin{equation*}
\left|\int^{t}_{0}\!\!\!\int_{\R} H(u_r)\,\rho_n'(x)\,dxds\right|\leq \|H\|_{\infty}t \int_{\R} |\rho_n'(x)|\,dx\le2\,\|H\|_{\infty}t\,,
\end{equation*}
we obtain
\begin{equation}\label{w2}
\left|u_s(\cdot,t)\left(\{x_j\}\right)- u_{0s}\left(\{x_j\}\right)\right| \leq 2\,\|H\|_{\infty}t\,.
\end{equation}
Since, by Theorem \ref{recl}, $u_s(\cdot,t)=u_s(\cdot,t)\lefthalfcup F$ and $F$ contains $p$ points, we obtain that $\|u_s(\cdot,t)-u_{0s}\|\le 2p\|H\|_\infty t$. Hence 
$u_s\in C([0,T];\mathcal{M}(\Omega))$ and $[u(\cdot,0)]_s=u_{0s}$ (see also Remark~\ref{continuity}).
 
If $u_{0s}\lefthalfcup \{x_j\}=\pm u_{0s}^\pm\lefthalfcup \{x_j\}$, it follows from \eqref{mono.us} that $[u(\cdot,t) ]_s^\mp\lefthalfcup \{x_j\}=0$ for any $t\in[0,T]$. Then inequality \eqref{w2} gives
$$
[u_s(\cdot,t)]^\pm\left(\{x_j\}\right)\geq u_{0s}^\pm\left(\{x_j\}\right)- 2t \|H\|_{\infty}>0$$
for  any $t\in[0,t_j)$, with $t_j:=\frac{u_{0s}^\pm\left(\{x_j\}\right)}{2\|H\|_{\infty}}$. 
Then by the monotonicity of the mappings $t\mapsto u_s^{\pm}(\cdot,t)$ (see \eqref{mono.us}) the conclusion  follows.
\hfill$\square$


\section{Problem $(D)$: comparison and uniqueness}\label{unid}
\setcounter{equation}{0}

As already said, to address $(CL)$ we need results concerning {\em singular} 
Dirichlet  initial-boundary value problems for the scalar conservation law:
\begin{equation*}
\left\{\begin{array}{ll}
u_{t}+[H(u)]_x=0&\mbox{in}\ \,\O\times(0,T)=:Q \\
u=m_1&\mbox{in $\{a\}\times (0,T)$} \\
u=m_2&\mbox{in $\{b\}\times (0,T)$} \\
u=u_0&\mbox{in $\Omega\times\{0\}$}\,,
\end{array}\right.\leqno{(D)}
\end{equation*}
where $\Omega=(a,b)$ is a bounded interval,  $m_1=\pm\infty$, $m_2=\pm\infty$, and $u_0:\O\mapsto\R$.    
Similar problems will be considered also for half-lines, either $\Omega\equiv(a,\infty)$, or $\Omega\equiv(-\infty,b)$; obviously, the above condition at $\{b\}\times (0,T)$ is omitted when $\Omega\equiv(a,\infty)$, and that at $\{a\}\times (0,T)$ is omitted when $\Omega\equiv(-\infty,b)$.

We shall denote problem $(D)$ by $(D_S)$ when $m_1=\pm\infty$, $m_2=\pm\infty$, or by $(D_R)$ when both $m_1$ and $m_2$ are finite. When $\Omega=(a,b)$ problem $(D_S)$ stands for four different initial-boundary value problems, which we denote by $(D_+^+)$, $(D_-^-)$, $(D_+^-)$ and $(D_-^+)$ according to the four choices $m_1=m_2=\infty$, $m_1=m_2=-\infty$, $m_1=\infty, m_2=-\infty$ and $m_1=-\infty, m_2=\infty$. In the case of half-lines problem $(D_S)$ consists only of two cases, namely
\begin{equation*}
\left\{\begin{array}{ll}
u_{t}+[H(u)]_x=0&\mbox{in}\ \,Q \\
u=\pm\infty&\mbox{in $\{a\}\times (0,T)$} \\
u=u_0&\mbox{in $\Omega\times\{0\}$}\,
\end{array}\right.\leqno{(D_\pm)}
\end{equation*}
if $\O=(a,\infty)$, and
\begin{equation*}
\left\{\begin{array}{ll}
u_{t}+[H(u)]_x=0&\mbox{in}\ \,Q \\
u=\pm\infty&\mbox{in $\{b\}\times (0,T)$} \\
u=u_0&\mbox{in $\Omega\times\{0\}$}
\end{array}\right.\leqno{(D^\pm)}
\end{equation*}
if $\O=(-\infty,b)$. We shall write that a statement holds  for problem $(D_S)$, if it collectively holds for all problems $(D_\pm^\pm)$. 
\smallskip

The following definition concerns problem $(D_R)$ (see \cite{Te}).

\begin{definition}\label{defre}
Let $\Omega=(a,b)$, $u_0\in BV(\Omega)$. 

\noindent $(i)$ An {\em entropy subsolution} of $(D_R)$ is any $\underline u\in  BV(Q)$ such that:

\noindent $(a)$  for every $k\in\R$ and for all $\zeta \in C^1([0,T];C^1_c(\Omega))$, $\zeta(\cdot,T)=0$ in $\Omega$, $\zeta\geq 0$ in $Q$, 
\begin{equation}\label{cond1+}
\iint_{Q}\!\left \{[ u-k]_+\zeta_t+{\rm sgn}_+( u-k)[H( u)-H(k)]\zeta_x\right\}dxdt \geq -\int_{\Omega}[u_0-k]_+\zeta(x,0)\,dx;
\end{equation}

\noindent $(b)$ for a.e. $t\in(0,T)$ there holds
\begin{equation}\label{sublimBV1--}
{\sgn}_+(\underline{u}(a^+,t)-k)
\left[H(\underline{u}(a^+,t))-H(k)\right]  \le 0\;\; \;\text{if  $k>m_1$}\,,
\end{equation}
\begin{equation}\label{sub2limBV--}
{\sgn}_+(\underline{u}(b^-,t)-k)\left[H(\underline{u}(b^-,t))-H(k)\right]\ge0 \;\; \;\text{if  $k>m_2$}\,.
\end{equation}

\noindent $(ii)$ An {\em entropy supersolution} of $(D_R)$ is any $\overline u\in BV(Q)$ such that:

\noindent $(a')$  for every $k\in\R$ and if  $\zeta \in C^1([0,T];C^1_c(\Omega))$, $\zeta(\cdot,T)=0$ in $\Omega$, $\zeta\geq 0$ in $Q$, 
\begin{equation}\label{cond1-}
\iint_{Q}\!\left \{[ u-k]_-\zeta_t+{\rm sgn}_-( u-k)[H( u)-H(k)]\zeta_x\right\}dxdt \geq -\int_{\Omega}[u_0-k]_-
\zeta(x,0)\,dx;
\end{equation}

\noindent $(b')$ for a.e. $t\in(0,T)$ there holds
\begin{equation}\label{superlim1BV++}
{\sgn}_-(\overline{u}(a^+,t)-k)\left[H(\overline{u}(a^+,t))-H(k)\right]  \le0
\;\;\;\text{if  $k<m_1$}\,,
\end{equation}
\begin{equation}\label{superlim2BV++}
{\sgn}_-(\overline{u}(b^-,t)-k)\left[H(\overline{u}(b^-,t))-H(k)\right] \ge0 \;\;\;\text{if  $k<m_2$}\,.
\end{equation}

\noindent $(iii)$ A function $\overline u\in BV(Q)\cap C([0,T];L^1(\Omega))$ is an {\em entropy solution} of $(D_R)$ if it is both an entropy subsolution and an entropy supersolution.
\end{definition}
When $\O=(a,\infty)$ entropy sub- and supersolutions of $(D_R)$ are defined as above, only dropping conditions \eqref{sub2limBV--} and \eqref{superlim2BV++}; similarly, conditions \eqref{sublimBV1--} and \eqref{superlim1BV++} are omitted if $\O=(-\infty,b)$. Moreover, in these cases we require that $\underline{u}, \overline{u}$ belong to $BV_{loc}(Q)\cap L^\infty(Q)$.

\begin{remark}\label{re01}
If  $\underline{u}, \overline{u} \in BV(Q)$, the traces $\underline{u}(a^+,t):={\rm ess} \lim_{\xi\to a^+} \underline{u}(\xi,t)$, $\underline{u}(b^-,t):={\rm ess} \lim_{\eta\to b^-} \underline{u}(\eta,t)$ exist for a.e. $t\in(0,T)$, and similarly for $\overline{u}$. Hence the above definitions are well-posed. By the same token, conditions \eqref{sublimBV1--}-\eqref{sub2limBV--} and \eqref{superlim1BV++}-\eqref{superlim2BV++} can be reformulated as follows: for every $\beta\in C^1_c(0,T)$, $\beta\ge0$, 
\begin{subequations}
\begin{equation}\label{wbv1}
{\rm ess} \lim_{\xi\to a^+} \int_0^T\!{\sgn}_+(\underline{u}(\xi,t)-k)\, \big[H(\underline{u}(\xi,t))-H(k)\big]\,\beta(t)\,dt  \le 0\; \;\text{if  $k>m_1$}\,,
\end{equation}
\begin{equation}\label{wbv2}
{\rm ess} \lim_{\eta\to b^-} \int_0^T\! {\sgn}_+(\underline{u}(\eta,t)-k)\, \big[H(\underline{u}(\eta,t))-H(k)\big]\,\beta(t)\,dt\ge0 \; \;\text{if  $k>m_2$}\,,
\end{equation}
\begin{equation}\label{wbv3}
{\rm ess} \lim_{\xi\to a^+} \int_0^T\!{\sgn}_-(\overline{u}(\xi,t)-k)\, \big[H(\overline{u}(\xi,t))-H(k)\big]\,\beta(t)\,dt  \le0\; \;\text{if  $k<m_1$}\,,
\end{equation}
\begin{equation}\label{wbv4}
{\rm ess} \lim_{\eta\to b^-}  \int_0^T\!{\sgn}_-(\overline{u}(\eta,t)-k)\, \big[H(\overline{u}(\eta,t))-H(k)\big]\,\beta(t)\,dt \ge0\; \;\text{if  $k<m_2$} \,.
\end{equation}
\end{subequations}
\end{remark}
The following definitions for problem $(D_S)$ are formulated for a wider class of initial data.
\begin{definition}\label{defsub}
Let $\O=(a,b)$, $u_0\in L^1(\Omega)$. 

\noindent $(i)$ An {\em entropy subsolution} of $(D_+^+)$ is any $\underline u\in C([0,T];L^1(\Omega))$ such that:

\noindent $(a)$
 for every $k\in\R$ and $\zeta \in C^1_c(Q)$,  $\zeta\geq 0$ in $Q$ 
\begin{equation}\label{sir}
\iint_{Q}\!\left \{[\underline u-k]_+\zeta_t\,+\,{\rm sgn}_+(\underline u-k)\,[H(\underline u)-H(k)]\zeta_x\right\}\,dxdt \geq 0\,, 
\end{equation}
and for any interval $I\subseteq\O$
\begin{equation}\label{sir0}
\lim_{t\to 0^+} \int_I \, [\underline u(x,t)-u_0(x)]_+\,dx=0 \,.
\end{equation}

\noindent $(ii)$ An {\em entropy subsolution} of $(D_-^-)$ is any $\underline u\in C([0,T];L^1(\Omega))$ such that
$(a)$ holds, and for every $k\in\R$, $\beta\in C^1_c(0,T)$, $\beta\ge0$, 
\begin{subequations}
\begin{equation}\label{sublim1--}
{\rm ess} \lim_{\xi\to a^+} \int_0^T\!{\sgn}_+(\underline{u}(\xi,t)-k)\, \big[H(\underline{u}(\xi,t))-H(k)\big]\,\beta(t)\,dt  \le 0\,,
\end{equation}
\begin{equation}\label{sub2lim--}
{\rm ess} \lim_{\eta\to b^-} \int_0^T\!{\sgn}_+(\underline{u}(\eta,t)-k)\, \big[H(\underline{u}(\eta,t))-H(k)\big]\,\beta(t)\,dt\ge0 \,.
\end{equation}
\end{subequations}
\noindent $(iii)$ An {\em entropy subsolution} of $(D_+^-)$ is any $\underline u\in C([0,T];L^1(\Omega))$ such that
$(a)$ holds, and for every $k\in\R$, $\beta\in C^1_c(0,T)$, $\beta\ge0$ inequality \eqref{sub2lim--} holds. 

\noindent $(iv)$ An {\em entropy subsolution} of $(D_-^+)$ is any $\underline u\in C([0,T];L^1(\Omega))$ such that
$(a)$ holds, and for every $k\in\R$, $\beta\in C^1_c(0,T)$, $\beta\ge0$ inequality \eqref{sublim1--} holds. 
\end{definition}
\begin{definition}\label{defsuper}
Let $\O=(a,b)$, $u_0\in L^1(\Omega)$. 

\noindent $(i)$ An {\em entropy supersolution} of $(D_+^+)$ is any $\overline u\in C([0,T];L^1(\Omega))$ such that: 

\noindent $(a')$ for every $k\in\R$ and $\zeta \in C^1_c(Q)$, $\zeta\geq 0$ in $Q$ 
\begin{equation}\label{sirbis}
\iint_{Q}\!\left \{[\overline u-k]_-\zeta_t+{\rm sgn}_-(\overline u-k)\,[H(\overline u)-H(k)]\zeta_x\right\}\,dxdt \geq 0\,, \end{equation}
and for any interval $I\subseteq\O$
\begin{equation}\label{sir0bis}
\lim_{t\to 0^+} \int_I \, [\overline u(x,t)-u_0(x)]_-\,dx=0 \,;
\end{equation}

\noindent $(b')$ for every $k\in\R$ and $\beta\in C^1_c(0,T)$, $\beta\ge0$, 
\begin{subequations}
\begin{equation}\label{superlim1++}
{\rm ess} \lim_{\xi\to a^+} \int_0^T\!{\sgn}_-(\overline{u}(\xi,t)-k)\, \big[H(\overline{u}(\xi,t))-H(k)\big]\,\beta(t)\,dt  \le0\,,
\end{equation}
\begin{equation}\label{superlim2++}
{\rm ess} \lim_{\eta\to b^-}  \int_0^T\!{\sgn}_-(\overline{u}(\eta,t)-k)\, \big[H(\overline{u}(\eta,t))-H(k)\big]\,\beta(t)\,dt \ge0 \,.
\end{equation}
\end{subequations}
\noindent $(ii)$ An {\em entropy supersolution} of $(D_-^-)$ is any $\overline u\in C([0,T];L^1(\Omega))$ such that
$(a')$ holds.

\noindent $(iii)$ An {\em entropy supersolution} of $(D_+^-)$ is any $\overline u\in C([0,T];L^1(\Omega))$ such that
$(a')$ holds, and for every $k\in\R$, $\beta\in C^1_c(0,T)$, $\beta\ge0$,  inequality \eqref{superlim1++} holds. 

\noindent $(iv)$ An {\em entropy supersolution} of $(D_-^+)$ is any $\overline u\in C([0,T];L^1(\Omega))$ such that
$(a')$ holds, and for every $k\in\R$, $\beta\in C^1_c(0,T)$, $\beta\ge0$,  inequality \eqref{superlim2++} holds. 
\end{definition}

\begin{definition}\label{defsol}
A function $u\in C([0,T];L^1(\Omega))$ is called an {\em entropy solution} of $(D_S)$ if it is both an entropy subsolution and an entropy supersolution of $(D_S)$. 
\end{definition}

Observe that \eqref{superlim1++}-\eqref{superlim2++} can be regarded as limiting cases of \eqref{wbv3}-\eqref{wbv4}, since for every $k\in\R$ there holds ${\rm sgn}_-(m_i-k)\to 0$ as $m_i\to\infty$ $(i=1,2)$. Similarly, \eqref{sublim1--}-\eqref{sub2lim--} can be regarded as limiting cases of \eqref{wbv1}-\eqref{wbv2} as $m_i\to-\infty$. 

\begin{remark}\label{rem.fd}
Let us prove that every entropy solution of $(D_S)$ satisfies the weak formulation 
\begin{equation}\label{fd.D}
\iint_Q \left\{u\zeta_t+H(u)\zeta_x\right\}\,dxdt=0
\end{equation}
for every $\zeta\in C^1_c(Q)$. 

To this aim, we fix any sequence $k_j\to -\infty$. By \eqref{sir}, for all $\zeta\in C^1_c(Q)$, $\zeta\geq 0$, there holds
\begin{equation}\label{fd1}
\iint_Q\left\{[u-k_j ]_+\zeta_t+ {\rm sgn}_+(u-k_j)\, [H(u)-H(k_j) ]\zeta_x\right\}\,dxdt \geq 0\,.
\end{equation}
Let us take the limit as $j\to \infty$ in \eqref{fd1}. Since $u\in L^1(Q)$ and $H$ is bounded, we have
\begin{eqnarray}\label{fd2}
&&\int_Q{\rm sgn}_+(u-k_j)\, [H(u)-H(k_j) ]\zeta_x\,dxdt= \iint_{\{u>k_j\}} H(u)\zeta_x\,dxdt\\
&&-\underbrace{\iint_{Q}H(k_j)\,\zeta_x\,dxdt}_{=0} +\iint_{\{u\leq k_j\}}H(k_j)\zeta_x\,dxdt\to  \iint_{Q} H(u)\zeta_x\,dxdt\,,\nonumber
\end{eqnarray}
and 
\begin{eqnarray}\label{fd3}
&&\iint_Q[u-k_j ]_+\zeta_t\,dxdt= \iint_{\{u>k_j\}}u\zeta_t\,dxdt -\overbrace{\iint_Q k_j\zeta_t\,dxdt}^{=0}\\
&& +\iint_{\{u\leq k_j\}}k_j\zeta_t\,dxdt \to \iint_Qu\zeta_t\,dxdt\,,\nonumber
\end{eqnarray}
(here we have used that $\left|\iint_{\{u\leq k_j\}}k_j\zeta_t\,dxdt \right|\leq \iint_{\{u\leq k_j\}}|u|\,|\zeta_t|\,dxdt \to 0$, as $k_j\to -\infty$). In view of \eqref{fd2}--\eqref{fd3}, letting $j\to \infty$ in \eqref{fd1} gives 
\begin{equation}\label{fd.D+}
\iint_Q \left\{u\zeta_t+H(u)\zeta_x\right\}\,dxdt\geq0
\end{equation}
for every $\zeta\in C^1_c(Q)$, $\zeta\geq 0$. Analogously, letting $k_j\to \infty$ in 
\begin{equation}\label{fd1bis}
\iint_Q\left\{[u-k_j ]_-\zeta_t+ {\rm sgn}_-(u-k_j)\, [H(u)-H(k_j) ]\zeta_x\right\}\,dxdt \geq 0\,
\end{equation} 
(see \eqref{sirbis}) gives, for every $\zeta$ as above, 
\begin{equation}\label{fd.D-}
\iint_Q \left\{u\zeta_t+H(u)\zeta_x\right\}\,dxdt\leq0\,.
\end{equation}
Therefore the conclusion follows combining \eqref{fd.D+} and \eqref{fd.D-}.
\end{remark}

\begin{remark} \label{comp-entr} 
The conditions (\ref{sublim1--}-\ref{sub2lim--}) and (\ref{superlim1++}-\ref{superlim2++}) are entropy boundary 
conditions for singular Dirichlet problems and give a meaning, 
in a hyperbolic sense, to the boundary conditions "$u=-\infty$" and "$u=\infty$".
As already mentioned in the Introduction, 
they coincide with the compatibility conditions \eqref{comp1} and \eqref{comp2} 
for entropy solutions of (CL) at points $x_j$ where a signed Dirac mass is concentrated. 
\end{remark}

\begin{remark}\label{conneHJ}
Let $u$ denote either $\underline{u}$ in \eqref{sir}, or $ \overline{u}$ in \eqref{sirbis}.  
Choosing $\zeta(x,t)=\alpha(x)\beta(t)$ with $\alpha\in C^1_c(\Omega)$, $\beta\in C^1_c(0,T)$,  $\alpha,\beta \ge0$, gives 
\begin{equation}\label{dim}
\iint_{Q} \! \! \left\{[u(x,t)\! -\! k]_{\pm}\,\alpha(x)\beta'(t)\! +\! {\rm sgn}_{\pm}(u(x,t)\! -\! k)\! \left [H(u(x,t))\! -\! H(k)\right ]\!  \alpha'(x)\beta(t)\right\}\! dxdt\!  \ge \! 0
\end{equation}
for any $k\in\R$ . Since $0\le [u_r-k]_{\pm}\le [u]_{\pm}+|k|$, from the above inequality we get 
\begin{eqnarray*}
&&-\! \int_{\Omega} dx\,\a'(x) \left\{\int_0^T \! \! {\rm sgn}_{\pm}(u(x,t)\! -\! k)\! \left [H(u(x,t))\! -\! H(k)\right ]\! \beta(t)\,dt\right\} \,  \le\\
&&\qquad \le   \|\beta'\|_\infty  \int_{\Omega}dx\,\a(x) \left\{\int_0^T \big([u]_{\pm}(x,t)+|k|\big)\,dt\right\}\,=\\
&&\qquad= -\,\|\beta'\|_\infty\! \!  \int_{\Omega}dx\,\a'(x) \left\{\int_0^T\!\!\! \int_{c}^x\big([u]_{\pm}(x,t)+|k|\big)\,dt\right\} \,.
\end{eqnarray*}
for every $c\in \overline{\Omega}$. Hence the  distributional derivative of the function 
\begin{equation*}
x\mapsto  \!\int_0^T\!\! {\rm sgn}_{\pm}(u(x,t)-k)\left [H(u(x,t))-H(k)\right]\beta(t)\,dt \,-\,  \|\beta'\|_\infty \!\int_0^T\!\!\!\int_{c}^x\big([u]_{\pm}(x,t)+|k|\big)\,dydt
\end{equation*}
is nonpositive. Therefore, the limits
\begin{subequations}
\begin{equation}\label{defili}
{\rm ess} \lim_{x\to a^+} \int_0^T \!\! {\rm sgn}_{\pm}(u(x,t)-k)\left [H(u(x,t))-H(k)\right]\beta(t)\,dt \,,
\end{equation}
\begin{equation}\label{defilibis}
{\rm ess} \lim_{x\to b^-} \int_0^T \!\! {\rm sgn}_{\pm}(u(x,t)-k)\left [H(u(x,t))-H(k)\right]\beta(t)\,dt \,
\end{equation}
\end{subequations}
exist and are finite, thus the above definitions are well-posed. 

The same statement can be applied to entropy solutions 
of $(CL)$, since they satisfy inequalities \eqref{dim} in every domain $Q_j=(x_j,x_{j+1})\times (0,T)$ ($j=1,\dots,p-1$), or $Q^-=(-\infty, x_1)\times (0,T)$, $Q^+=(x_p,\infty)\times (0,T)$ (recall that by  Theorem \ref{recl} and assumption $(A_1)$ the singular part of an entropy solution of $(CL)$ is not supported in these domains).
\end{remark}
\begin{remark}\label{re02}
Conditions \eqref{sublim1--}-\eqref{sub2lim--} for subsolutions of $(D_-^-)$ can be equivalently 
rewritten as follows: for all $k\in\R$ and $\b$ as above and for a.e.~$\xi,\eta\in(a,b)$
\begin{subequations}
\begin{equation}\label{sub1-}
\int_0^T\!\!\!\!\int_a^\xi [\underline{u}(x,t)-k]_+\beta'(t)\,dxdt \,   
\ge  \int_0^T\!{\sgn}_+(\underline{u}(\xi,t)-k)\, \big[H(\underline{u}(\xi,t))-H(k)\big]\,\beta(t)\,dt  \end{equation}
\begin{equation}\label{sub2-}
\int_0^T\!\!\!\!\!\int_\eta^b [\underline{u}(x,t)-k]_+\beta'(t)\,dxdt 
\ge - \!\!\int_0^T\!{\sgn}_+(\underline{u}(\eta,t)-k)\left[H(\underline{u}(\eta,t))-H(k)\right]\beta(t)\,dt.
\end{equation}
\end{subequations}
Similarly, conditions \eqref{superlim1++}-\eqref{superlim2++} for supersolutions of $(D_+^+)$ equivalently 
read: for all $k\in\R$ and for a.e.~$\xi,\eta \in(a,b)$
\begin{subequations}
\begin{equation}\label{super1++}
\int_0^T\!\!\!\!\!\int_a^\xi [\overline{u}(x,t)-k]_-\beta'(t)\,dxdt  \ge 
\int_0^T\!{\sgn}_-(\overline{u}(\xi,t)-k) \left[H(\overline{u}(\xi,t))-H(k)\right]\,\beta(t)\,dt 
\end{equation}
\begin{equation}\label{super2++}
\int_0^T\!\!\!\!\!\int_\eta^b [\overline{u}(x,t)-k]_-\beta'(t)\,dxdt  \ge  
- \int_0^T\!{\sgn}_-(\overline{u}(\eta,t)-k) \left[H(\overline{u}(\eta,t))-H(k)\right]\beta(t)\,dt .
\end{equation}
\end{subequations}
\end{remark}

When $\O=(a,\infty)$ we have the following definition (we omit the formulation for the case $\O=(-\infty,b)$).
\begin{definition}\label{-infb}
Let $\O=(a,\infty)$, $u_0\in L^1_{loc}(\Omega)$. 

\noindent $(i)$ An {\em entropy subsolution} of $(D_+)$ is any $\underline u\in C([0,T];L^1_{loc}(\Omega))$ such that $(a)$ of Definition \ref{defsub} holds.

\noindent $(ii)$ An {\em entropy subsolution} of $(D_-)$ is any $\underline u\in C([0,T];L^1_{loc}(\Omega))$ such that $(a)$ of Definition \ref{defsub} holds, and for every $k\in\R$, $\beta\in C^1_c(0,T)$, $\beta\ge0$ inequality \eqref{sublim1--} holds.

\noindent $(iii)$ An {\em entropy supersolution} of $(D_+)$ is any $\overline u\in C([0,T];L^1_{loc}(\Omega))$ such that $(a')$ of Definition \ref{defsuper} holds, and for every $k\in\R$, $\beta\in C^1_c(0,T)$, $\beta\ge0$,  inequality \eqref{superlim1++} holds. 

\noindent $(iv)$ An {\em entropy supersolution} of $(D_-)$ is any $\overline u\in C([0,T];L^1_{loc}(\Omega))$ such that
$(a')$ of Definition \ref{defsuper} holds. 

\noindent $(v)$ A function $u\in C([0,T];L^1_{loc}(\Omega))$ is called an {\em entropy solution} of $(D_S)$ if it is both an entropy subsolution and an entropy supersolution of $(D_S)$. 
\end{definition} 

Comparison and uniqueness results for problem $(D_R)$ are given by the following theorem (see \cite[Theorem 1.1]{Te}). 
\begin{theorem}\label{cf1} 
Let $\O=(a,b)$. Let $u_0,v_0\in BV(\Omega)$, and $m_1,m_2,n_1,n_2\in\R$. Let $\underline u$ be an entropy subsolution of $(D_R)$, and $\overline v$  be an entropy supersolution of $(D_R)$ with $u_{0}$, $m_1$ and $m_2$ replaced by $v_{0}$, $n_1$ and $n_2$. Then for a.e. $t\in(0,T)$ 
\begin{equation}\label{co1}
\int_{\Omega}[\underline u(x,t)-\overline v(x,t)]_+\,dx \le \int_{\Omega}[u_0(x)- v_0(x)]_+\,dx \,+\, 
 \big([m_1-n_1]_+ + [m_2-n_2]_+ \big)\|H'\|_{\infty}\,t\; .
\end{equation}
Similar  results hold for $\O=(a,\infty)$ and $\O=(-\infty,b)$ if $u_0,v_0\in BV_{loc}(\Omega)\cap L^\infty(\O)$. 
In these cases for a.e. $t\in(0,T)$ there holds
\begin{equation}\label{co01}
\int_a^R[\underline u(x,t)-\overline v(x,t)]_+\,dx \le \int_a^{R+\|H'\|_{\infty}t}[u_0(x)- v_0(x)]_+\,dx \,+\, [m_1-n_1]_+ \|H'\|_{\infty}\,t
\end{equation}
 for every $R>a$ if $\O=(a,\infty)$, respectively
  \begin{equation*}
\int^b_R[\underline u(x,t)-\overline v(x,t)]_+\,dx \le \int^b_{R-\|H'\|_{\infty}t}[u_0(x)- v_0(x)]_+\,dx \,+\, 
 [m_2-n_2]_+ \|H'\|_{\infty}\,t
\end{equation*}
 for every $R<b$ if  $\O=(-\infty,b)$. Therefore, in all cases there exists at most one solution of $(D_R)$.
\end{theorem}

As for problem $(D_S)$, the following holds.
 \begin{theorem}\label{cmp} Let $(A_0)$ hold. Let $\underline{u},\, \overline{u}$ be an entropy sub- and supersolution of $(D_S)$ with the same boundary conditions. Then $\underline{u}\le\overline{u}$ a.e. in $Q$. In particular, there exists at most one entropy solution of $(D_S)$.
\end{theorem}
\begin{proof} We only give the proof for $(D_+^-)$, as in the other cases of $(D_S)$ it is similar. 
We use the Kru\v zkov doubling method adapted to boundary valued problems (see \cite{BSTT2, MNRR, O, Se}). Let $\rho_{\epsilon}$ $(\epsilon>0)$ be a symmetric mollifier in $\R$, and set 
\begin{equation*}
\z(x,t,y,s):=\,\rho_{\epsilon_1}(x-y)\,\rho_{\epsilon_2}(t-s)\, \s_1\Big(\frac{x+y}2\Big) \,\s_2\Big(\frac{t+s}2\Big)\qquad ((x,t),\,(y,s)\in Q)\,,
\end{equation*} 
with $\s_1\in C^1_c(\O)$, $\s_2\in C^1_c(0,T)$, $\s_1 \ge0$, $\s_2 \ge0$. From  \eqref{sirbis} and \eqref{sir} we get
\begin{eqnarray*}
&&\iint_{Q} \big\{\sgn_-(\overline{u}(x,t)-\underline{u}(y,s))[H(\overline{u}(x,t))-H(\underline{u}(y,s))]\z_x(x,t,y,s)\,+ \\
&&\qquad \quad+[ \overline{u}(x,t)-\underline{u}(y,s)]_-\,\z_t(x,t,y,s)\big\}\, dxdt \, \ge 0\quad\text{for all $(y,s)\in Q$} \,, \nonumber 
\end{eqnarray*}
\begin{eqnarray*}
&&\iint_{Q} \big\{\sgn_+(\underline{u}(y,s)-\overline{u}(x,t))[H(\underline{u}(y,s))-H(\overline{u}(x,t))]\z_y(x,t,y,s) \,+  \\
&&\qquad \quad+[\underline{u}(y,s)\!-\! \overline{u}(x,t)]_+\,\z_s(x,t,y,s)\big\}\,dyds\ge  0
\quad\text{for all $(x,t)\in Q$}\,. \nonumber 
\end{eqnarray*} 
Recalling that $[u]_+=[-u]_-$ and $\sgn_+(u)=-\sgn_-(-u)$ ($u\in\R$), we sum  the above inequalities
integrated over $Q$\,:
\begin{equation}\label{xtys}
\begin{aligned}
&\iint\!\!\!\!\!\!\iint_{Q\times Q} \rho_{\epsilon_1}(x-y)\,\rho_{\epsilon_2}(t-s)
\left\{[\overline{u}(x,t)-\underline{u}(y,s)]_- \, \s_1\left(\frac{x+y}2\Big) \,\s_2'\Big(\frac{t+s}2\right)+\right. \\
&\left.+  \sgn_-(\overline{u}(x,t)\!-\! \underline{u}(y,s))[H(\overline{u}(x,t))\!-\! H(\underline{u}(y,s))]   \s_1'\left(\frac{x\!+\!y}2\right) \s_2\left(\frac{t\!+\!s}2\right)\right\} dxdtdyds \ge 0. 
\end{aligned}
\end{equation}
Set
\begin{eqnarray*}
I_1:=\iint\!\!\!\!\!\!\iint_{Q\times Q} \rho_{\epsilon_1}(x-y)\,\rho_{\epsilon_2}(t-s)
[\overline{u}(x,t)-\underline{u}(y,s)]_- \, \s_1(y) \,\s_2'\Big(\frac{t+s}2\Big) \,dxdtdyds\,,
\end{eqnarray*}
$$
\begin{aligned}
&I_2:=\iint\!\!\!\!\!\!\iint_{Q\times Q} \rho_{\epsilon_1}(x-y)\,\rho_{\epsilon_2}(t-s) \,\times \\
&\qquad\times\;\sgn_-(\overline{u}(x,t)- \underline{u}(y,s))[H(\overline{u}(x,t))- H(\underline{u}(y,s))] \,  \s_1'(y)\,\s_2\Big(\frac{t+s}2\Big)\,dxdtdyds\, .
\end{aligned}
$$
Observe that the difference between $I_1$ and the first term in \eqref{xtys} vanishes as $\epsilon_1\to 0^+$; the same holds for the difference between $I_2$ and the second term in \eqref{xtys}.

Let $a<\xi<\eta<b$ be fixed. By standard approximation arguments we can choose 
$$
\s_1(y)\equiv \s_{1,n}(y)=n(y-\xi)\chi_{[\xi,\xi+1/n]}(y)
+\chi_{(\xi+1/n,\eta-1/n)}(y)-n(y-\eta)\chi_{[\eta-1/n,\eta]}(y),
$$
where $n\in\mathbb{N}$ and $ y\in\O$,
thus
$$
\s_1'(y)=n\chi_{[\xi,\xi+1/n]}(y)-n\chi_{[\eta-1/n,\eta]}(y)\,.
$$
With this choice of $\s_1$, $I_2$ reads
\begin{eqnarray}\label{I21}
&& \qquad  I_2
=\,n\int_0^T\!\!\!ds\!\int_\xi^{\xi+1/n} \!\! dy \!\int_a^b \!\! dx \,\rho_{\epsilon_1}(x-y) \, \times \\
&&\times \int_0^T \!\!dt \;\sgn_-(\overline{u}(x,t) - \underline{u}(y,s))[H(\overline{u}(x,t)) -H(\underline{u}(y,s))] \,\rho_{\epsilon_2}(t-s)\,\s_2\Big(\frac{t+s}2\Big)\, - \nonumber \\
&&-\,n\int_0^T\!\!\!ds\!\int_{\eta-1/n}^\eta \!\! dy\!\int_a^b \!\! dx \,\rho_{\epsilon_1}(x-y) \, \times \nonumber \\
&&\times \int_0^T \!\!dt \;\sgn_-(\overline{u}(x,t) - \underline{u}(y,s))[H(\overline{u}(x,t)) -H(\underline{u}(y,s))] \,\rho_{\epsilon_2}(t-s)\,\s_2\Big(\frac{t+s}2\Big)\, . \nonumber 
\end{eqnarray}
By \eqref{super1++} and \eqref{sub2-}, from \eqref{I21} we obtain for a.e. $\xi,\eta\in(a,b)$
\begin{eqnarray*}
&&
\lim_{n\to\infty}I_2= \int_0^T\!\!\!ds \!\int_a^b \!\! dx \,\rho_{\epsilon_1}(x-\xi) \, \times  \\
&&\times \int_0^T \!\!dt \;\sgn_-(\overline{u}(x,t) - \underline{u}(\xi,s))[H(\overline{u}(x,t)) -H(\underline{u}(\xi,s))] \,\rho_{\epsilon_2}(t-s)\,\s_2\Big(\frac{t+s}2\Big) \, - \nonumber \\
&&-\int_0^T\!\!ds\!\!\int_a^b \!\! dx \,\rho_{\epsilon_1}(x-\eta) \, \times \nonumber \\
&&\times \int_0^T \!\!dt \;\sgn_-(\overline{u}(x,t) - \underline{u}(\eta,s))[H(\overline{u}(x,t)) -H(\underline{u}(\eta,s))] \,\rho_{\epsilon_2}(t-s)\,\s_2\Big(\frac{t+s}2\Big)\, \le \nonumber \\
&& \le \int_0^T\!\!\!ds \!\int_a^b \!\! dx \,\rho_{\epsilon_1}(x-\xi) \int_0^T \!\!dt \! \int_a^x \!\!dz \,
[\overline{u}(z,t)-\underline{u}(\xi,s)]_- \left[\rho_{\epsilon_2}(t-s)\,\s_2\Big(\frac{t+s}2\Big)\right]_t \,+\nonumber \\
&&+ \int_0^T\!\!\!dt \!\int_a^b \!\! dx \,\rho_{\epsilon_1}(x-\eta) \int_0^T \!\!ds \! \int_\eta^b \!\!dz \,
[\underline{u}(z,s)-\overline{u}(x,t)]_+ \left[\rho_{\epsilon_2}(t-s)\,\s_2\Big(\frac{t+s}2\Big)\right]_s \,=:S \,. \nonumber
\end{eqnarray*}
As  $\epsilon_1\to 0^+$ we obtain that 
\begin{eqnarray}\label{I23}
&& \lim_{\epsilon_1\to 0^+}S =
\int_0^T\!\!\!ds \!\int_0^T \!\!dt \! \int_a^\xi \!\!dz \,
[\overline{u}(z,t)-\underline{u}(\xi,s)]_- \left[\rho_{\epsilon_2}(t-s)\,\s_2\Big(\frac{t+s}2\Big)\right]_t \,+ \\
&&\qquad + \int_0^T\!\!\!dt \!\int_0^T \!\!ds \! \int_\eta^b \!\!dz \,
[\underline{u}(z,s)-\overline{u}(\eta,t)]_+ \left[\rho_{\epsilon_2}(t-s)\,\s_2\Big(\frac{t+s}2\Big)\right]_s \,\le \nonumber\\
&& \qquad\le C_{\epsilon_2} T\, \Big\{ \int_0^T \!\!dt \! \int_a^\xi\!\!dz\, |\overline{u}(z,t)|  \,+\,(\xi-a)\! \!\int_0^T |\underline{u}(\xi,s)|\,ds  \,+ \nonumber \\
&& \qquad+ \int_0^T \!\!ds \! \int_\eta^b\!\! dz\, |\underline{u}(z,s)|  \,+\,(b-\eta)\! \!\int_0^T |\overline{u}(\eta,t)|\,dt  \Big\} \,. \nonumber
\end{eqnarray}
Clearly, there holds
\begin{equation*}
\lim_{\xi\to a^+}\int_0^T \!\!dt \! \int_a^\xi\!\!dz\, |\overline{u}(z,t)| \, =
\lim_{\eta\to b^-}\int_0^T \!\!ds \! \int_\eta^b\!\!dz\, |\underline{u}(z,s)| =0\,.
\end{equation*}
On the other hand, since $\int_0^T |\underline{u}(\xi,s)|\,ds\ge0$ and the map $\xi \to \int_0^T |\underline{u}(\xi,s)|\,ds$ belongs to $L^1(\O)$, there holds
$$
{\rm ess} \liminf\limits_{\xi\to a^+} \;(\xi-a)\! \!\int_0^T |\underline{u}(\xi,s)|\,ds=0
$$
- for, otherwise there would exist $c,\d>0$ such that $(\xi-a)\! \!\int_0^T |\underline{u}(\xi,s)|\,ds\ge c$ for a.e. $\xi\in(a,a+\d)$, thus $\xi\to\int_0^T |\underline{u}(\xi,s)|\,ds\not \in L^1(\O)$. Therefore, for every $n\in\N$ there exist $\d_n>0$ and $E_n\subseteq(a,a+\d_n)$,  $|E_n|>0$,  such that $(\xi-a)\!\int_0^T |\underline{u}(\xi,s)|\,ds<1/n$ for a.e. $\xi\in E_n$. It follows that a sequence $\{\xi_n\}\subseteq\O$ exists, such that:

\smallskip

\noindent $(i)$ $\xi_n$ is a Lebesgue point of $\int_0^T |\underline{u}(\xi,s)|\,ds$,

\smallskip

\noindent $(ii)$ $\xi_n\to a^+$ as $n\to\infty$, and $(\xi_n-a)\!\int_0^T |\underline{u}(\xi_n,s)|\,ds<\frac1n$ for all $n\in\N$\,.

Similarly, there holds
$$
{\rm ess} \liminf\limits_{\eta\to b^-} \;(b-\eta)\! \!\int_0^T |\underline{u}(\eta,t)|\,dt=0\,,
$$
hence there exists $\{\eta_n\}\subseteq\O$, $\eta_n\to b^-$ as $n\to\infty$, with properties analogous to $(i)$-$(ii)$ above. Then writing \eqref{I23} with $\xi=\xi_n$, $\eta=\eta_n$ and letting  $n\to\infty$ we obtain that the right-hand side of  \eqref{I23} goes to zero.

To sum up, following the above procedure and letting $\e_2\to0^+$ from \eqref{xtys}, we get for any $\s_2\in C^1_c(0,T)$, $\s_2 \ge0$,
\begin{equation}\label{I24}
\iint_{Q} [\overline{u}(x,t)-\underline{u}(x,t)]_- \, \s_2'(t)\,dxdt\ge0\,.
\end{equation}
Let $0<t_1<t_2\leq T$ be fixed. By standard approximation arguments we can choose 
$$
\s_2(t)\equiv \s_{2,n}(t)=n(t-t_1)\chi_{[t_1,t_1+1/n]}+\chi_{(t_1+1/n,t_2-1/n)}-n(t-t_2)\chi_{[t_2-1/n,t_2]}\qquad (n\in\mathbb{N})\,.
$$
Then from \eqref{I24} we get for all $n$
$$
n\int_{t_2-1/n}^{t_2}\int_{\Omega} [\overline{u}(x,t)-\underline{u}(x,t)]_-\,dxdt \, \le \, n \int_{t_1}^{t_1+1/n}\int_{\Omega} [\overline{u}(x,t)-\underline{u}(x,t)]_-\,dt\,,
$$
whence as $n\to \infty$ 
$$
\int_{\Omega}  [\overline{u}(x,t_2)-\underline{u}(x,t_2)]_-\,dx \leq \int_{\Omega}  [\overline{u}(x,t_1)-\underline{u}(x,t_1)]_-\,dx\,.
$$
Since $\underline{u},\overline{u}\in C([0,T];L^1(\O))$, as $t_1\to 0^+$\ by \eqref{sir0} and \eqref{sir0bis} there holds $$
\int_{\Omega}  [\overline{u}(x,t_2)-\underline{u}(x,t_2)]_-\,dx = 0 \quad \text{for all $t_2\in(0,T]$}\,.
$$
This proves the result.
\end{proof}
 
For future reference we prove the following generalization of \cite[Lemma 4.4]{BSTT2}.
\begin{prop}\label{letra}
\noindent $(i)$ Let $u$ be an entropy solution either of $(D_+^+)$, or  of $(D_+^-)$. Then there exists $f_a^+\in L^{\infty}(0,T)$ such that
\begin{equation}\label{z1}
{\rm ess}\lim_{x\to a^+}\int_0^T H(u(x,t))\,\beta(t)\,dt=\int_0^T f_a^+(t)\,\beta(t)\,dt\,
\end{equation}
for every $\beta\in C^1_c(0,T)$, and 
\begin{equation}\label{x1}
\limsup_{u\to \infty} H(u)\leq f_a^+(t)\leq \sup_{u\in\R}H(u)\quad\mbox{for}\ \,a.e.\ \,t\in (0,T)\,.
\end{equation}

\noindent $(ii)$ Let $u$ be an entropy solution either of $(D_-^+)$, or  of $(D_-^-)$. Then there exists $f_a^-\in L^{\infty}(0,T)$ such that
\begin{equation}\label{z3}
{\rm ess}\lim_{x\to a^+} \int_0^T H(u(x,t))\,\beta(t)\,dt=\int_0^T f_a^-(t)\,\beta(t)\,dt\,
\end{equation}
for every $\beta\in C^1_c(0,T)$, and
\begin{equation}\label{x3}
\inf_{u\in\R} H(u)\leq f_a^-(t)\leq \liminf_{u\to -\infty} H(u)\quad\mbox{for}\ \,a.e.\ \,t\in (0,T)\,.
\end{equation}

\noindent $(iii)$ Let $u$ be an entropy solution either of $(D_-^+)$, or  of $(D_+^+)$. Then there exists $f_b^+\in L^{\infty}(0,T)$ such that
\begin{equation}\label{z2}
{\rm ess}\lim_{x\to b^-} \int_0^T H(u(x,t))\,\beta(t)\,dt=\int_0^T f_b^+(t)\,\beta(t)\,dt\,
\end{equation}
for every $\beta\in C^1_c(0,T)$, and 
\begin{equation}\label{x2}
\inf_{u\in\R} H(u)\leq f_b^+(t)\leq \liminf_{u\to \infty} H(u)\quad\mbox{for}\ \,a.e.\ \,t\in (0,T)\,.
\end{equation}

\noindent $(iv)$ Let $u$ be an entropy solution either of $(D_+^-)$, or  of $(D_-^-)$. Then there exists $f_b^-\in L^{\infty}(0,T)$ such that
\begin{equation}\label{z4}
{\rm ess}\lim_{x\to b^-} \int_0^T H(u(x,t))\,\beta(t)\,dt=\int_0^T f_b^-(t)\,\beta(t)\,dt\,
\end{equation}
for every $\beta\in C^1_c(0,T)$, and 
\begin{equation}\label{x4}
\limsup_{u\to -\infty} H(u)\leq f_b^-(t)\leq \sup_{u\in\R}H(u)\quad\mbox{for}\ \,a.e.\ \,t\in (0,T)\,.
\end{equation}
\end{prop}

\begin{proof}
The existence of the limits in the left-hand side of \eqref{z1}, \eqref{z3}, \eqref{z2} and \eqref{z4} follows from \eqref{defili}-\eqref{defilibis}, since for $a.e.$ $x\in\Omega$ there holds  
\begin{eqnarray*}
x\mapsto \int_0^T H(u(x,t))\beta(t)\,dt = \int_0^T [{\rm sgn}_+(u(x,t)) -{\rm sgn}_-(u(x,t))]H(u(x,t))\beta(t)\,dt\, \end{eqnarray*}
(recall that $H(0)=0$). On the other hand, for every sequence $\{x_n\}$, $x_n\to a^+$, the sequence $\{H(u(x_n,\cdot))\}$ is bounded in $L^{\infty}(0,T)$. Hence there exist a subsequence $\{x_{n_k}\}\subseteq \{x_n\}$ and a function $f_a^+\in L^{\infty}(0,T)$ (independent of $\{x_{n_k}\}$) such that $H(u(x_{n_k},\cdot))\stackrel{*}\rightharpoonup f_a^+$ in $L^{\infty}(0,T)$, thus \eqref{z1} follows. Equalities \eqref{z3}, \eqref{z2} and \eqref{z4} are similarly proven.

Let us prove \eqref{x1}. Clearly, there holds $f_a^+(t)\le \sup_{u\in\R} H(u)$ for a.e. $t\in (0,T)$.
To prove the first inequality, let us choose in \eqref{sir} $\zeta(x,t)=\rho(x)\beta(t)$ with $\rho \in C^1_c([a,b))$, $\rho\geq 0$, $\beta\in C^1_c(0,T)$, $\beta\geq 0$. By standard arguments we can also choose $\rho=\a\s_{\epsilon}$ with $\a \in C^1_c([a,b))$, $\a\geq 0$, and
\begin{equation}\label{defeta}
\s_{\epsilon}(x):=\frac{2(x-a)-\epsilon}{\epsilon}\chi_{[a+\epsilon/2,a+\epsilon]}(x)+\chi_{(a+\epsilon,b]}(x)\qquad (x\in \Omega)\,.
\end{equation}
Then for every $k\in \R$ we obtain that
\begin{eqnarray*}
&&\iint_{Q}\Big\{ [u(x,t)-k]_+\a(x)\s_{\epsilon}(x)\b'(t)\,+\\
&&\qquad + {\rm sgn}_+(u(x,t)- k)[H(u(x,t))-H(k)] \,\a'(x)\s_{\epsilon}(x)\b(t)\Big\}\,dxdt\ge\\
&&\quad \ge\; -\frac2\epsilon\int_0^T\!\!dt\,\b(t)\!\int_{a+\epsilon/2}^{a+\epsilon} {\rm sgn}_+(u(x,t)- k)[H(u(x,t))-H(k)]\,\a(x)\,dx\,.
\end{eqnarray*}
Letting $\epsilon\to 0^+$ and using \eqref{superlim1++} and \eqref{z1}, we get that for every $k\in \R$ 
\begin{eqnarray*}
&& \iint_{Q}\left\{ [u(x,t)\!-\!k]_+\a(x)\b'(t)+ {\rm sgn}_+(u(x,t)\!-\! k)[H(u(x,t))\!-\!H(k)] \a'(x)\b(t)\right\}dxdt \ge \\
&&\qquad \ge -\a(a)\,{\rm ess}\lim_{x\to a^+}\int_0^T\!\! {\rm sgn}_+(u(x,t)-k)[H(u(x,t))-H(k)]\beta(t)\,dt\,= \\
&&\qquad =-\a(a)\,\Big\{ {\rm ess}\lim_{x\to a^+}\int_0^T [H(u(x,t))-H(k)]\,\beta(t)\,dt\,+\\
&&\qquad+ \underbrace{{\rm ess}\lim_{x\to a^+}\int_0^T\!\! {\rm sgn}_-(u(x,t)-k)[H(u(x,t))-H(k)]\,\beta(t)\,dt}_{\leq 0}\Big\}\,\ge \\
&&\qquad \ge -\a(a) \int_0^T [f_a^+(t)-H(k)]\,\beta(t)\,dt\,. 
\end{eqnarray*}
Letting $k\to \infty$ in the above inequality gives
\begin{equation*}
0\le \int_0^T [f_a^+(t)-\limsup_{k\to \infty} H(k)]\,\beta(t)\,dt\,,
\end{equation*}
whence by the arbitrariness of $\b$ inequality \eqref{x1} follows. 

To prove \eqref{x3} we argue as for \eqref{x1}, using inequality \eqref{sublim1--},  \eqref{sirbis} and \eqref{z3} instead of \eqref{sir}, \eqref{superlim1++} and \eqref{z1}. Then we get for every $k\in \R$ 
\begin{eqnarray*}
&& \iint_{Q}\left\{ [u(x,t)\!-\!k]_-\a(x)\b'(t)\!+\! {\rm sgn}_-(u(x,t)\!-\! k)[H(u(x,t))\!-\!H(k)] \,\a'(x)\b(t)\right\}dxdt \ge \\
&&\qquad\ge  -\a(a)\,{\rm ess}\lim_{x\to a^+}\int_0^T\!\! {\rm sgn}_-(u(x,t)\!-\!k)[H(u(x,t))\!-\!H(k)]\beta(t)\,dt= \\
&&\qquad=-\a(a)\,\Big\{- {\rm ess}\lim_{x\to a^+}\int_0^T [H(u(x,t))-H(k)]\,\beta(t)\,dt\,+\\
&&\qquad+ \underbrace{{\rm ess}\lim_{x\to a^+}\int_0^T\!\! {\rm sgn}_+(u(x,t)-k)[H(u(x,t))-H(k)]\,\beta(t)\,dt}_{\le 0}\Big\}\,\ge \\
&&\qquad\ge  \a(a) \int_0^T [f_a^-(t)-H(k)]\,\beta(t)\,dt\,. 
\end{eqnarray*}
As $k\to -\infty$ in the above inequality, by the arbitrariness of $\b$ we obtain
$$
f_a^-(t)\leq \liminf_{k\to -\infty} H(k)\quad\mbox{for}\ \,a.e.\ \,t\in (0,T)\,,
$$
thus \eqref{x3} follows.
 The proof of \eqref{x2} and \eqref{x4} is similar to that  of \eqref{x1} and \eqref{x3}, using
\begin{equation}\label{defetabis}
\s_{\epsilon}(x):=\chi_{[a,b-\epsilon)}(x)-\frac{2(x-b)+\epsilon}{\epsilon}\chi_{[b-\epsilon,b-\epsilon/2]}(x)\qquad (x\in \Omega)\,
\end{equation}
instead of \eqref{defeta}; we leave the details to the reader.
\end{proof}
Finally we prove the following result.
\begin{lem}\label{esti} Let $u$ be an entropy solution of $(D_R)$. Then for every $t\in (0,T]$
\begin{equation}\label{el1}
\|u(\cdot,t)\|_{L^1(\Omega)}\leq \|u_0\|_{L^1(\Omega)}+2\,\|H\|_{\infty}t\,. 
\end{equation}
\end{lem} 
\begin{proof} By \eqref{cond1+} and \eqref{cond1-} there holds 
\begin{equation*}
\iint_{Q} \left\{|u-k|\zeta_s+{\rm sgn}(u-k)\,[H(u)-H(k)]\zeta_x\right\}\,dxds \geq -\int_{\Omega}|u_0-k|\,\zeta(x,0)\,dx
\end{equation*}
for every $k\in\R$ and $\zeta$ as above. By standard arguments we can choose $\zeta(x,s)=\a_p(x)\b_q(s)$ with
$$
\a_p(x)=p(x-a)\chi_{[a,a+1/p)}(x)+\chi_{[a+1/p,b-1/p)}(x)-p(x-b)\chi_{[b-1/p, b]}(x)\,,
$$
$$
\b_q(s)=\chi_{[0,t-1/q)}(s)-q(s-t)\chi_{[t-1/q,t)}(s)\,
$$
 for any fixed $t\in(0,T]$ and $p,q\in\N$ sufficiently large.
Then for $k=0$ as $q\to\infty$  we get
$$\int_{\Omega}|u(x,t)|\,\a_p(x)\,dx\,-\int_{\Omega}|u_0(x)|\,\a_p(x)\,dx 
\le 2\,\|H\|_{\infty}t\,, 
$$
whence as $p\to\infty$ \eqref{el1} follows.
\end{proof}


\section{Problem $(D)$: existence}\label{exid}
\setcounter{equation}{0}

Let us recall the following result (see \cite{BLN, Te}).
\begin{theorem}\label{euBV}
Let $\Omega=(a,b)$, and let $u_0\in BV(\Omega)$. Then there exists a unique entropy solution $u\in BV(Q)\cap C([0,T];L^1(\Omega))$ of problem $(D_R)$. Moreover,
\begin{equation}\label{ep11}
\|u\|_{L^{\infty}(Q)}\leq \max\left\{|m_1|,|m_2|,\|u_0\|_{L^{\infty}(\Omega)}\right\}\,.
\end{equation}

\noindent The same holds for\, $\O=(a,\infty)$ and $\O=(-\infty,b)$ with  $u_0\in BV_{loc}(\Omega)\cap L^\infty(\O)$, ${\rm supp} \, u_0$ compact. In these cases there holds $u\in BV_{loc}(Q)\cap L^\infty(Q)\cap C([0,T];L^1_{loc}(\Omega))$, and inequality \eqref{ep11} is replaced by
\begin{equation}\label{ep11a}
\|u\|_{L^{\infty}(Q)}\leq \max\left\{|m_1|,\|u_0\|_{L^{\infty}(\Omega)}\right\}
\end{equation}
if $\O=(a,\infty)$, respectively by $\|u\|_{L^{\infty}(Q)}\leq \max\left\{|m_2|,\|u_0\|_{L^{\infty}(\Omega)}\right\}$ if $\O=(-\infty,b)$.
\end{theorem}
The above uniqueness claim follows from Theorem \ref{cf1}. Let us outline the proof of the existence part; we limit ourselves to the case $\O=(a,\infty)$, the proof being the same for $\O=(-\infty,b)$ and easier for $\Omega=(a,b)$.

Let $f_{1,\ep},f_{2,\ep} \in C^{\infty}(\R)$ ($0<\ep<1$) be a partition of unity:
\begin{equation*}
\left\{\begin{array}{ll}f_{1,\ep}=1\quad\mbox{in}\ \,(-\infty,a+2\sqrt{\ep} ]\,,\quad {\rm supp}\,f_{1,\ep}\subseteq (-\infty, a+3\sqrt{\ep} ]\smallskip\\
f_{2,\ep}=1\quad\mbox{in}\ \,[a+3\sqrt{\ep},\infty)\,,\quad {\rm supp}\,f_{2,\ep}\subseteq[a+2\sqrt{\ep}, \infty)\,,\smallskip\\
0\leq f_{i,\ep}\leq 1\,,\;\; \sum_{i=1}^2 f_{i,\ep}=1 \quad\mbox{in}\ \,\R\,,
\end{array}\right.
\end{equation*}
such that, for  $i=1,2$,
\begin{equation*}
\sup_{\ep\in (0,1)}\|f'_{i,\ep}\|_{L^1(\R)}<\infty\,, \quad\sup_{\ep \in (0,1)} \sqrt{\ep}\,\|f'_{i,\ep}\|_{{L^{\infty}(\R)}}<\infty\,,\quad
\sup_{\ep \in (0,1)}\sqrt{\ep} \|f''_{i,\ep}\|_{L^1(\R)}<\infty.
\end{equation*}
Let $u_0\in BV_{loc}(\Omega)\cap L^\infty(\O)$ have compact support. Set
\begin{equation}\label{au0}
u_{0\ep }:=m_1 f_{1,\ep}+[\s_{\ep}*u_0 ]f_{2,\ep}\,, 
\end{equation}
where $\{\s_{\ep}\}$ is a family of standard mollifiers with supp$\,\s_{\ep}\subseteq [-\sqrt{\ep},\sqrt{\ep}]$. Then there holds  $u_{0\ep}\in C^{\infty}(\R)$, $u_{0\ep}=m_1$ in $[a, a+\sqrt{\ep}]$, ${\rm supp} \, u_{0\ep}$ compact. Moreover, 
\begin{equation}\label{st0}
\sup_{\ep \in (0,1)}\|u_{0\ep}\|_{L^{\infty}(\Omega)}\leq \max\left\{|m_1|,\|u_0\|_{L^{\infty}(\Omega)}\right\}\,,
\end{equation}
\begin{equation}\label{bv1}
\sup_{\ep \in (0,1)}
\| u_{0\ep}'\|_{L^1(\Omega)}<\infty \,,\,\quad 
\sup_{\ep \in (0,1)} \ep\, \| u_{0\ep}''\|_{L^1(\Omega)}<\infty
\end{equation}
\begin{equation}\label{st1}
 u_{0\ep}\to u_0\quad\mbox{in $L^p(\Omega)$}\quad \mbox{for every $p\in[1,\infty)$}\,, \qquad
 u_{0\ep}\stackrel{*}\rightharpoonup u_0\quad\mbox{in $L^{\infty}(\Omega)$}\,.
\end{equation}

Let $H\in W^{1,\infty}(\R)$, $H(0)=0$. Set 
\begin{equation*}
H_{\ep}(u):= g_{\ep}(u)\, \big([\s_{\ep}*H ](u)-[\s_{\ep}*H ](0)\big) \qquad\ \ (u\in\R)\,,
\end{equation*}
where the family $\{ g_{\ep}\}\in C^{\infty}_c(\R)$ satisfies $ g_{\ep}=1$ in $(-1/\ep ,1/\ep)$, $0\leq g_{\ep}(x)\leq 1$ in $\R$, supp$\, g_{\ep}\subseteq (-2/\ep,2/\ep)$. It is easily seen that 
\begin{equation}\label{coHe}
\left\{
\begin{array}{ll}
H_{\ep}(0)=0\,,\quad \|H_{\ep}\|_{W^{1,\infty}(\R)}\leq \|H\|_{W^{1,\infty}(\R)}\,, \smallskip
\\ 
\text{$H_{\ep}\to H$ \quad\text{uniformly on the compact subsets of $\R$}\,.
}\end{array}
\right. 
\end{equation}

 Let $u_{\ep}\in C^{2,1}(\overline{Q})$ be the unique  classical solution of the parabolic problem 
 \begin{equation*}
\left\{\begin{array}{ll}
u_{\ep t}+[H_{\ep}(u_{\ep}) ]_x=\ep u_{\ep xx}&\mbox{in}\ \,Q\smallskip\\
u_{\ep}=m_1&\mbox{in $\{a\}\times (0,T)$}\,\smallskip\\
u_{\ep}=u_{0\ep}&\mbox{in $\Omega\times\{0\}$}\,,
\end{array}\right.\leqno{(D_\ep)}
\end{equation*}
with $m_1\in\R$, $u_{0\ep}$ and $H_{\ep}$ as above 
($e.g.$, see \cite{LSU}). 
\begin{lem} There holds
\begin{equation}\label{ep1}
\sup_{\ep \in (0,1)}\|u_{\ep}\|_{L^{\infty}(Q)}\leq 
\max\left\{|m_1|,\|u_0\|_{L^{\infty}(\Omega)}\right\}\,,
\end{equation}
and there exists $c>0$ only depending on $m_1$, $TV(u_0;\Omega)$, and $\|H\|_{W^{1,\infty}(\R)}$) such that
\begin{equation}\label{ep0}
\sup_{\ep \in (0,1)}\|u_{\ep x}\|_{L^{\infty}(0,T;L^1(\Omega))}\,\le\, c \,,
\end{equation}
\begin{equation}\label{ep2}
\sup_{\ep \in (0,1)}\|u_{\ep t}\|_{L^{\infty}(0,T;L^1(\Omega))}\,\le\, c\,,
\end{equation}
\begin{equation}\label{ep3}
\sup_{\ep \in (0,1)}\left(\ep\|u_{\ep x}\|_{L^{\infty}(Q)}\right)\,\le\, c\,.
\end{equation}
\end{lem}
\begin{proof}
Inequality \eqref{ep1} follows by the maximum principle and \eqref{st0}.
Arguing as in the proof of \cite[Proposition 3.1]{Te}  (see also \cite{BLN}) and using \eqref{bv1} gives \eqref{ep0}-\eqref{ep2}. As for \eqref{ep3}, integrating the first equation of $(D_\ep)$ over $(a,x)$ gives
\begin{equation}\label{ta1}
\ep u_{\ep x}(x,t)-\ep u_{\ep x}(a,t)= \int_a^x\!\!u_{\ep t}(y,t)\,dy +H_{\ep}(u_{\ep}(x,t)) -H_{\ep}(m_1)\,,
 \end{equation}
whence
$$
\ep\, |u_{\ep x}(a,t)|\le \int_a^x\!\!|u_{\ep t}(y,t)|\,dy +2\|H\|_{\infty}+  |u_{\ep x}(x,t)|\,.
$$
Integrating the above inequality over $(a,a+1)$ and using \eqref{ep0}-\eqref{ep2} we get
\begin{equation}\label{ta2}
\ep\, |u_{\ep x}(a,t)|\le 2\|H\|_{\infty}+\tilde{c}
 \end{equation}
for some $\tilde{c}>0$ independent of $\ep$.
 Then by \eqref{ep0}-\eqref{ep2} and \eqref{ta1}-\eqref{ta2} the estimate in \eqref{ep3} follows.
\end{proof}

\noindent {\em Proof of Theorem \ref{euBV}.}
By  estimates \eqref{ep1}-\eqref{ep2} the family $\{u_{\ep}\}$ is bounded in $L^{\infty}(Q)$, and there exists $M>0$ (only depending on $m_1$, $TV(u_0;\Omega)$, $\|H\|_{W^{1,\infty}(\R)}$) such that
$$
\sup_{\ep \in (0,1)} \|u_{\ep x}\|_{L^{\infty}(0,T;L^1(\Omega))}+\sup_{\ep \in (0,1)}\|u_{\ep t}\|_{L^{\infty}(0,T;L^1(\Omega))}\leq M\,.
$$
Then by embedding theorems there exist a sequence $\{u_{\ep_n}\}\subseteq\{u_\ep\}$ and a function $u\in BV(Q)\cap C([0,T];L^1(\Omega))$ such that
\begin{equation}\label{ccl1}
u_{\ep_n}\to u\quad\mbox{in}\ \,C([0,T];L^1(\Omega)) \quad\mbox{as}\ \,n\to \infty\,.
\end{equation}
Arguing as in \cite{BLN} shows that $u$ is an entropy solution of problem $(D_R)$, In fact, let $E,F_{\ep}:\R\to\R$,  $E\in C^2(\R)$, $F_{\ep}\in C^1(\R)$ and $F'_{\ep}=E'H'_{\ep}$\,. Multiplying the first equation in $(D_\ep)$ by $E'(u_{\ep})\,\zeta$ gives for any $\zeta \in C^1([0,T];C^1_c(\overline{\O}))$, $\z(\cdot,T)=0$ in $\O$,
\begin{eqnarray}\label{dipa}
&& \quad-\int_{\Omega}E(u_{0\ep})(x)\zeta(x,0)\,dx +\ep\iint_{Q}\left\{E''(u_{\ep}) u_{\ep x}^2\zeta\,+\,[E(u_{\ep}) ]_x \zeta_x\right\}\,dxdt=\\
&=& \iint_{Q} \left\{E(u_{\ep})\zeta_t+F_{\ep}(u_{\ep})\zeta_x\right\}dxdt + 
\int_0^T \{F_{\ep}(m_1)-\ep E'(m_1)u_{\ep x}(a,t)\}\zeta(a,t)\,dt \,. \nonumber 
\end{eqnarray} 
By standard regularization arguments we can choose in \eqref{dipa} $E(u_{\ep})=[u_{\ep}-k]_\pm$, thus obtaining for all $k\in\R$ and $\z$ as above, $\z\ge0$,
\begin{eqnarray}\label{dipabis}
&&\iint_{Q}\!\left \{[u_{\ep}-k]_\pm\zeta_t+{\rm sgn}_\pm(u_{\ep}-k)\,[H_\ep(u_{\ep})-H_\ep(k)]\zeta_x\right\}\,dxdt\, \ge \\
&&\qquad \ge\ep \iint_{Q} {\rm sgn}_\pm(u_{\ep}-k)\, u_{\ep x}\,\zeta_x\,dxdt  -\int_{\Omega}[u_{0\ep}-k]_{ \pm}\,\zeta(x,0)\,dx \,- \nonumber\\
&&\qquad- {\rm sgn}_\pm(m_1-k) \int_0^T [H_{\ep}(m_1)-H_{\ep}(k)-\ep u_{\ep x}(a,t)]\zeta(a,t)\,dt \,. \nonumber  
\end{eqnarray} 
If $\z(\cdot,t)\in C^1_c(\Omega)$ for all $t\in(0,T)$, from \eqref{dipabis} we obtain 
\begin{eqnarray}\label{dipa3}
&&\iint_{Q}\!\left \{[u_{\ep}-k]_\pm\zeta_t+{\rm sgn}_\pm(u_{\ep}-k)\,[H_\ep(u_{\ep})-H_\ep(k)]\zeta_x\right\}\,dxdt\, \ge \\
&&\qquad \ge\ep \iint_{Q} {\rm sgn}_\pm(u_{\ep}-k)\, u_{\ep x}\,\zeta_x\,dxdt  -\int_{\Omega}[u_{0\ep}-k]_\pm\,\zeta(x,0)\,dx \,. \nonumber
\end{eqnarray} 
On the other hand, choosing in \eqref{dipabis} $\zeta(x,t)=\chi_{[a,\,\xi+1/n)}(x)\beta(t)$ $(\xi\in\O, n\in\N)$ with $\beta\in C^1_c(0,T)$, $\beta\ge0$, and letting $n\to\infty$ plainly gives for every $k\in\R$ 
\begin{eqnarray}\label{dipater}
&&\int_0^T\!\!\!\!\int_a^\xi [u_\ep(x,t)-k]_\pm\beta'(t)\,dxdt \,- 
\\&&\qquad - 
\int_0^T\!{\sgn}_\pm(u_\ep(\xi,t)-k)\, \big[H_\ep(u_\ep(\xi,t))-H_\ep(k)\big]\,\beta(t)\,dt\, \ge  \nonumber\\
&&\qquad \ge\ep \int_0^T\!\! {\rm sgn}_\pm(u_{\ep}(\xi,t)-k)\, u_{\ep x}(\xi,t)\,\b(t)dt  \,- \nonumber\\
&&\qquad - \,{\rm sgn}_\pm(m_1-k) \int_0^T [H_{\ep}(m_1)-H_{\ep}(k)-\ep u_{\ep x}(a,t)]\b(t)\,dt \,.\nonumber 
\end{eqnarray} 
Multiplying the first equation of $(D_\ep)$ by $\zeta(x,t)=\chi_{[a,\,\xi+1/n)}(x)\beta(t)$  and letting $n\to\infty$,  one easily sees that 
\begin{eqnarray}\label{dipaqua}
&&-\,\ep\int_0^T u_{\ep x}(a,t)\b(t)\,dt  =
-\, \ep\int_0^T u_{\ep x}(\xi,t)\b(t)\,dt -\\
&&\qquad - \int_0^T\!\!\!\!\int_a^\xi [u_\ep(x,t)-k]\beta'(t)\,dxdt\, + 
\int_0^T\! \big[H_\ep(u_\ep(\xi,t))-H_\ep(m_1)\big]\,\beta(t)\,dt\,.  \nonumber 
\end{eqnarray} 
From \eqref{dipater}-\eqref{dipaqua} we get 
\begin{eqnarray}\label{dipaqui}
&&\int_0^T\!\!\!\!\int_a^\xi\left\{ [u_\ep(x,t)-k]_\pm-{\rm sgn}_\pm(m_1-k) [u_\ep(x,t)-k] \right\}\beta'(t)\,dxdt \,\ge \\
&&\quad \ge \int_0^T\!\!\left[{\sgn}_\pm(u_\ep(\xi,t)\!-\!k)\!-\!{\rm sgn}_\pm(m_1\!-\!k)\right] 
\left[H_\ep(u_\ep(\xi,t))\!-\!H_\ep(k)\right]\beta(t)\,dt+  \nonumber\\
&&\quad +\,\ep \int_0^T \left[{\rm sgn}_\pm(u_{\ep}(\xi,t)-k)
{\rm sgn}_\pm(m_1-k)\right]\, u_{\ep x}(\xi,t)\,\b(t)dt  \,. \nonumber 
\end{eqnarray}

By \eqref{coHe}, \eqref{ep0} and \eqref{ccl1} we can take the limit as $\ep_n\to0^+$ in \eqref{dipa3} and \eqref{dipaqui} (written with $\ep=\ep_n$). It follows that the function $u$ in
\eqref{ccl1} satisfies the following inequalities:
 
\noindent - for every $k\in\R$ and for all $\zeta \in C^1([0,T];C^1_c(\Omega))$, $\zeta(\cdot,T)=0$ in $\Omega$, $\zeta\geq 0$ in $Q$, 
\begin{equation*}
\iint_{Q}\!\left \{[ u-k]_\pm\zeta_t\,+\,{\rm sgn}_\pm( u-k)\,[H( u)-H(k)]\zeta_x\right\}\,dxdt \geq -\int_{\Omega}[u_0-k]_\pm\,\zeta(x,0)\,dx\,;
\end{equation*}

\noindent - for every $k\in\R$ and $\beta\in C^1_c(0,T)$, $\beta\ge0$ and for a.e. $\xi\in\O$,
\begin{eqnarray*}
&&\int_0^T\!\!\!\!\int_a^\xi\left\{ [u(x,t)-k]_\pm-{\rm sgn}_\pm(m_1-k) [u(x,t)-k] \right\}\beta'(t)\,dxdt \,\ge \\
&&\qquad \ge \int_0^T\!\left[{\sgn}_\pm(u(\xi,t)-k)-{\rm sgn}_\pm(m_1-k)\right]\, \big[H(u(\xi,t))-H(k)\big]\,\beta(t)\,dt\,. \nonumber 
\end{eqnarray*} 
Letting $\xi\to a^+$ in the latter inequality and using Remark \ref{re01} we conclude that $u$ is an entropy solution of $(D_R)$. Hence the result follows. 
\hfill$\square$
\begin{remark} 
In the proof of Theorem \ref{euBV} when $\Omega=(a,b)$ one uses the family of solutions of the problem
 \begin{equation*}
\left\{\begin{array}{ll}
u_{\ep t}+[H_{\ep}(u_{\ep}) ]_x=\ep u_{\ep xx}&\mbox{in}\ \,Q\smallskip\\
u_{\ep}=m_1&\mbox{in $\{a\}\times (0,T)$}\,\smallskip\\
u_{\ep}=m_2&\mbox{in $\{b\}\times (0,T)$}\, \smallskip\\
u_{\ep}=u_{0\ep}&\mbox{in $\Omega\times\{0\}$}\,,
\end{array}\right.\leqno{(D_\ep')}
\end{equation*}
with $m_1,m_2\in\R$, $H_{\ep}$ as above and $u_{0\ep}$ defined by a suitable partition of unity; we leave the details to the reader.
\end{remark}

Concerning $(D_S)$ the following holds.
 \begin{theorem}\label{ecs}  Let $(A_0)$ hold. When $\O=(a,b)$ for any $u_0\in L^1(\Omega)$ there exists an entropy solution of $(D_S)$.  The same holds for any $u_0\in L^1_{loc}(\Omega)$ if $\O=(a,\infty)$, or $\O=(-\infty,b)$.
\end{theorem}
\begin{proof} Let $\O=(a,b)$. Let us prove the result for $(D^-_+)$, the proof being the same for $(D^+_+)$, $(D^-_-)$ and $(D^+_-)$.
Let $u_0\in BV(\Omega)$. By Theorem \ref{euBV}, for all $n,p\in\N$ there exists an entropy solution $u_{n,-p}\in BV(Q)\cap C([0,T];L^1(\Omega))$ of problem $(D_R)$ with $m_1=n$, $m_2=-p$. In particular, there holds:

\smallskip

\noindent $(a)$ by \eqref{sir}-\eqref{sir0} and \eqref{sirbis}-\eqref{sir0bis}, for every $k\in\R$ and for all $\zeta \in C^1([0,T];C^1_c(\Omega))$, $\zeta(\cdot,T)=0$ in $\Omega$, $\zeta\geq 0$ in $Q$,
\begin{eqnarray}\label{tir1}
&&\iint_{Q}\!\left \{[u_{n,-p}-k]_\pm\zeta_t+
{\rm sgn}_\pm(u_{n,-p}-k)\,[H(u_{n,-p})-H(k)]\zeta_x\right\}\,dxdt \,\ge \\
&&\qquad \ge -\int_{\Omega}[u_0-k]_\pm\,\zeta(x,0)\,dx\,; \nonumber
\end{eqnarray}

\smallskip

\noindent $(b)$  by \eqref{super1++}, for every  $\beta\in C^1_c(0,T)$, $\beta\ge0$, 
for $a.e.$ $\xi\in \Omega$ and for all $k<n$ and $p\in\R$,
\begin{subequations}
\begin{eqnarray}\label{tedi2}
&&\int_0^T\!\!\!\!\int_a^\xi [u_{n,-p}(x,t)-k]_-\beta'(t)\,dxdt \, \ge \\
&&\qquad \ge \int_0^T\!{\sgn}_-(u_{n,-p}(\xi,t)-k)\, \big[H(u_{n,-p}(\xi,t))-H(k)\big]\,\beta(t)\,dt;
\nonumber
\end{eqnarray}
$(c)$  by \eqref{sub2-}, for every  $\beta\in C^1_c(0,T)$, $\beta\ge0$, for $a.e.$ $\eta\in \Omega$
and for all $n\in\R$ and $k>-p$,
\begin{eqnarray}\label{tedi3}
&&\int_0^T\!\!\!\!\int_\eta^b [u_{n,-p}(x,t)-k]_+\beta'(t)\,dxdt \, \ge \\
&&\qquad \ge -\int_0^T\!{\sgn}_+(u_{n,-p}(\eta,t)-k)\, \big[H(u_{n,-p}(\eta,t))-H(k)\big]\,\beta(t)\,dt;
\nonumber
\end{eqnarray}
\end{subequations}

\smallskip

\noindent $(d)$  by \eqref{el1}, for every $t\in (0,T)$
\begin{equation}\label{po1}
\|u_{n,-p}(\cdot,t)\|_{L^1(\Omega)}\leq \|u_0\|_{L^1(\Omega)}+2\,\|H\|_{\infty}t\,. 
\end{equation}
 Moreover, by inequality \eqref{co1}, for all $n,p\in\N$ there holds a.e. in $Q$
\begin{subequations}\label{pnem}
\begin{equation}\label{nem}
u_{n,-p}\leq u_{n+1,-p}\,, 
\end{equation}
\begin{equation}\label{pem}
u_{n,-p}\geq u_{n,-p-1}\,. 
\end{equation}
\end{subequations}

Let $p\in\R$ be fixed. By \eqref{po1} and \eqref{nem} there exists $u_{\infty,-p}  \in L^{\infty}(0,T;L^1(\Omega))$ such that
\begin{equation}\label{po2}
u_{n,-p}\to u_{\infty,-p}\quad\text{in $L^1(Q)$ as $n\to\infty$}.
\end{equation}
Then letting $n\to\infty$ in \eqref{tir1} gives
\begin{eqnarray}\label{tir10}
&&\iint_{Q}\!\left \{[u_{\infty,-p}-k]_\pm\zeta_t+
{\rm sgn}_\pm(u_{\infty,-p}-k)\,[H(u_{\infty,-p})-H(k)]\zeta_x\right\}\,dxdt \,\ge \\
&&\qquad \ge  -\int_{\Omega}[u_0-k]_\pm\,\zeta(x,0)\,dx\,, \nonumber
\end{eqnarray}
whereas from \eqref{tedi2}-\eqref{tedi3} we get, for every  $\beta\in C^1_c(0,T)$, $\beta\ge0$, 
for $a.e.$ $\xi,\eta\in \Omega$ and for all $k,p\in\R$:
\begin{subequations}
\begin{eqnarray}\label{tedi20}
&&\int_0^T\!\!\!\!\int_a^\xi [u_{\infty,-p}(x,t)-k]_-\beta'(t)\,dxdt \, \ge \\
&&\qquad \ge \int_0^T\!{\sgn}_-(u_{\infty,-p}(\xi,t)-k)\, \big[H(u_{\infty,-p}(\xi,t))-H(k)\big]\,\beta(t)\,dt\,,
\nonumber
\end{eqnarray}
\begin{eqnarray}\label{tedi30}
&&\int_0^T\!\!\!\!\int_\eta^b [u_{\infty,-p}(x,t)-k]_+\beta'(t)\,dxdt \, \ge \\
&&\quad \ge -\!\int_0^T\!\!{\sgn}_+(u_{\infty,-p}(\eta,t)\!-\!k) \left[H(u_{\infty,-p}(\eta,t))\!-\!H(k)\right]\beta(t)\,dt \quad \text{($k>-p$)}\,.
\nonumber
\end{eqnarray}
\end{subequations}

Moreover, from \eqref{po1} and \eqref{pem} we obtain
\begin{equation}\label{po10}
\|u_{\infty,-p}(\cdot,t)\|_{L^1(\Omega)}\leq \|u_0\|_{L^1(\Omega)}+\|H\|_{\infty}t\quad
\text{for every $t\in (0,T)$\,,} 
\end{equation}
\begin{equation}\label{pem0}
u_{\infty,-p}\geq u_{\infty,-p-1}\quad\text{a.e. in $Q$}\,. 
\end{equation}

By \eqref{po10}-\eqref{pem0} there exists $u_{\infty,-\infty}\in L^{\infty}(0,T;L^1(\Omega))$ such that
\begin{equation}\label{po20}
u_{\infty,-p}\to u_{\infty,-\infty}\quad\text{in $L^1(Q)$ 
as $p\to\infty$}.
\end{equation}
Then letting $p\to\infty$ in \eqref{tir10} shows that $u_{\infty,-\infty}$ satisfies \eqref{sir} and \eqref{sirbis}. In addition, letting $p\to\infty$ in \eqref{tedi20}-\eqref{tedi30} proves that $u_{\infty,-\infty}$ satisfies \eqref{sub2lim--} and \eqref{superlim1++} for every $k\in\R$ (see Remark \ref{re02}). 
By Remark \ref{continuity} and arguing as in the proof of \cite[Proposition 3.20]{BSTT1}, it can be checked that $u_{\infty,-\infty}\in C([0,T];L^1(\Omega))$ and, by construction, $u_{\infty,-\infty}(\cdot,0)=u_0$ in $\mathcal{M}(\Omega)$. Therefore \eqref{sir0} and \eqref{sir0bis} follow as well, and $u_{\infty,-\infty}$ is an entropy solution of $(D_+^-)$. 

\smallskip

It remains to remove the assumption $u_0\in BV(\O)$. To this purpose, let $v_0\in BV(\O)$, and let $v_{\infty,-\infty}$ be the entropy solution of  $(D_+^-)$ with initial data $v_0$ constructed by the above procedure  (and with the same boundary conditions considered in the construction of $u_{\infty,-\infty}$). Then by \eqref{co1} there holds
\begin{equation}\label{co11}
\|u_{\infty,-\infty}(\cdot,t)-v_{\infty,-\infty}(\cdot,t)\|_{L^1(\Omega)}\leq \|u_0-v_0\|_{L^1(\Omega)}\quad\mbox{for every}\ \,t\in (0,T)\,.
\end{equation}	
Let $u_0\in L^1(\Omega)$, let $\{u_{0j}\}\subseteq BV(\Omega)$ be any sequence such that
$u_{0j}\to u_0$ in $L^1(\Omega)$.
Let $\{u_j\}\equiv\{(u_{\infty,-\infty})_j\}$ be the sequence of entropy solutions to problem $(D_+^-)$ constructed as above, with initial data $u_{0j}$. Then for all $j\in\N$:

\smallskip

\noindent $(a)$ for every $k\in\R$ and for all $\zeta \in C^1([0,T];C^1_c(\Omega))$, $\zeta(\cdot,T)=0$ in $\Omega$, $\zeta\geq 0$ in $Q$,
\begin{eqnarray}\label{tirj}
&&\iint_{Q}\!\left \{[u_j-k]_\pm\zeta_t+
{\rm sgn}_\pm(u_j-k)\,[H(u_j)-H(k)]\zeta_x\right\}\,dxdt \,\ge \\
&&\qquad \ge -\int_{\Omega}[u_{0j}-k]_\pm\,\zeta(x,0)\,dx\,; \nonumber
\end{eqnarray}

\noindent $(b)$ for all $\beta\in C^1_c(0,T)$, $\beta\ge0$, for all $k\in\R$  and for $a.e.$ $\xi,\eta\in \Omega$:
\begin{subequations}
\begin{eqnarray}\label{tedj}
&&\int_0^T\!\!\!\!\int_a^\xi [u_j(x,t)-k]_-\beta'(t)\,dxdt \, \ge \\
&&\qquad \ge \int_0^T\!{\sgn}_-(u_j(\xi,t)-k)\, \big[H(u_j(\xi,t))-H(k)\big]\,\beta(t)\,dt\,,
\nonumber
\end{eqnarray}
\begin{eqnarray}\label{tedjbis}
&&\int_0^T\!\!\!\!\int_\eta^b [u_j(x,t)-k]_+\beta'(t)\,dxdt \, \ge \\
&&\qquad \ge -\int_0^T\!{\sgn}_+(u_j(\eta,t)-k)\, \big[H(u_j(\eta,t))-H(k)\big]\,\beta(t)\,dt \,.
\nonumber
\end{eqnarray}
\end{subequations}

By \eqref{co11} there holds 
\begin{equation*}
\|u_i-u_j\|_{L^1(Q)}\leq T\|u_{0i}-u_{0j}\|_{L^1(\Omega)}\quad\mbox{for all}\ \,i,j\in\mathbb{N}\,,
\end{equation*}	
thus there exists $u \in L^1(Q)$ such that $u_j\to u$ in $L^1(Q)$ as $j\to\infty$. As before, there holds $u\in C([0,T];L^1(\Omega))$.  Then letting $j\to\infty$ in \eqref{tirj} and \eqref{tedj}-\eqref{tedjbis} the result for $(D_+^-)$ follows. The other cases of $(D_S)$ can be dealt with similarly, hence the conclusion follows if $\O=(a,b)$.

The above arguments easily extend to the case of half-lines. For instance, let $\O=(a,\infty)$, $m_1=\infty$ and $u_0\in BV(\Omega)\cap L^\infty(\O)$. Then by Theorem \ref{euBV} and inequality \eqref{co01} there exists a sequence $\{u_n\}$ of entropy solutions of $(D_R)$, such that for every $t\in (0,T)$ $\|u_n(\cdot,t)\|_{L^1(\Omega)}\leq \|u_0\|_{L^1(\Omega)}+2\,\|H\|_{\infty}t$, and $u_n\leq u_{n+1}$ a.e. in $Q$ for all $n\in\N$\,. Then letting $n\to\infty$ we obtain an entropy solution $u_\infty$ of  $(D_+)$ in this case. Moreover, if $v_0\in BV(\Omega)\cap L^\infty(\O)$ and $v_\infty$ is the corresponding entropy solution of  $(D_+)$ with initial data $v_0$ constructed as before, by \eqref{co01} there holds
\begin{equation}\label{co011}
\|u_\infty(\cdot,t)-v_\infty(\cdot,t)\|_{L^1(a,R)}\leq \|u_0-v_0\|_{L^1(a,R+\|H'\|_\infty t)}\quad\mbox{for every}\ \,t\in (0,T)\,.
\end{equation}	

Now let $u_0\in L^1_{loc}(\Omega)$, and let $\{u_{0j}\}\subseteq BV(\Omega)\cap L^\infty(\O)$, ${\rm supp} \, u_{0j}$ compact, $u_{0j}\to u_0$ in $L^1_{loc}(\Omega)$.
Let $\{u_j\}\equiv\{(u_{\infty})_j\}$ be the sequence of entropy solutions of $(D_+)$ constructed as above, with initial data $u_{0j}$. By \eqref{co011} $\{u_j\}$ is a Cauchy sequence in $L^\infty(0,T;L^1(K))$ for every compact subset $K\subset\O$. Then by a diagonal argument the conclusion easily follows. 
\end{proof}


\section{Well-posedness of problem $(CL)$}  
\label{cleu}
\setcounter{equation}{0}

In this section we prove Theorem \ref{exiuni}. 
\medskip

We first prove the existence claim. Rewrite $(A_1)$ as follows:
\begin{equation}\label{H0r}
u_{0s}=\sum_{j=1}^r c_j^+\delta_{x_j'}-\sum_{j=1}^s c_j^-\delta_{x_j''}
\qquad 
( c_j^\pm>0, \,r+s= p)\,.
\end{equation}
For every $j=1,\dots,p$ such that $c_j>0$ we set
\begin{equation}\label{C_j+}
C_j^+(t):=\Big[\,  c_j - \int_0^t \left[f_{x_j^+}^+(s) - f_{x_j^-}^+(s) \right] ds
\,\Big]_+ \qquad (t\in [0,T])\,,
\end{equation}
with $f_{x_j^+}^+$ satisfying \eqref{z1} (written with $x_j^+$ instead of $a$) and $f_{x_j^-}^+$ satisfying \eqref{z2} (written with $x_j^-$ instead of $b$); observe that by \eqref{x1} and \eqref{x2} there holds 
\begin{equation}\label{f_j+}
f_{x_j^+}^+(s) - f_{x_j^-}^+(s)\ge0 \;\;\text{ for a.e. $s\in(0,T)$}\,.
\end{equation}
 Similarly, for every $j=1,\dots,p$ such that $c_j<0$ we set
\begin{equation}\label{C_j-}
C_j^-(t):=\Big[\,  c_j - \int_0^t \left[f_{x_j^+}^-(s) - f_{x_j^-}^-(s) \right] ds
\,\Big]_- \qquad (t\in [0,T])\,,
\end{equation}
with $f_{x_j^+}^-$ satisfying \eqref{z3} (written with $x_j^+$ instead of $a$) and $f_{x_j^-}^-$ satisfying \eqref{z4} (written with $x_j^-$ instead of $b$); observe that by \eqref{x3} and \eqref{x4} there holds 
\begin{equation}\label{f_j-}
f_{x_j^+}^-(s) - f_{x_j^-}^-(s)\le0\;\;\text{ for a.e. $s\in(0,T)$}\,.
\end{equation}

Let $\bar{t}_j:=\sup \{t\in [0,T]\,|\, C_j^\pm(t)>0\}>0$ $(j=1,\dots,p)$. Then $\bar{t}_j>0$ since $C_j^\pm(0)= \pm c_j>0$. By \eqref{f_j+} and \eqref{f_j-}  $C_j^\pm$ is nonincreasing in $(0,T)$, whence $C_j^\pm>0$ in $[0,\bar{t}_j)$ and, if $\bar{t}_j<T$, there holds $C_j^\pm=0$ in $[\bar{t}_j,T]$. Let $\tau_1:=\min\,\{\bar{t}_1,\dots,\bar{t}_p \}$, and define $u\in C([0,\tau_1]; \ \mathcal M(\R))$ as follows: 
\begin{subequations}\label{defuu}
\begin{equation}\label{defur}
\begin{cases}
\text{ in $Q_{1,\t_1}$ $u_r$ is the entropy solution  of $(D^+)$ if $c_1>0$, of $(D^-)$ if $c_1<0$;} \\
\text{ in $Q_{j,\t_1}$ ($j=2,\dots,p$) $u_r$ is the entropy solution  of $(D_+^+)$ if $\min\{c_{j-1},c_j\}>0$,} \\
\text{ of $(D_-^-)$ if $\max\{c_{j-1},c_j\}<0$, of $(D_+^-)$ if $c_{j-1}>0>c_j$,  of $(D_-^+)$ if $c_{j-1}<0<c_j$;} \\
\text{ in $Q_{p+1,\t_1}$ $u_r$ is the entropy solution  of $(D_+)$ if $c_p>0$, of $(D_-)$ if $c_p<0$} \,.
\end{cases}
\end{equation}
\begin{equation}\label{defus}
u_s(\cdot,t):=\sum_{j=1}^r C_j^+(t)\delta_{x_j'}-\sum_{j=1}^s C_j^-(t)\delta_{x_j''}\;.
\end{equation}
\end{subequations}
By Definitions \ref{enso} and \ref{defsub}-\ref{defsol} $u$ is an entropy solution of $(CL)$ in $Q_{j,\tau_1}$ for $j=1,\dots,p+1$. Hence $u$ is an entropy solution of $(CL)$ in $S_{\tau_1}$, if we prove \eqref{ewf}-\eqref{mkru} with $\Omega=\R$, $\tau=\tau_1$ for all $\zeta\in C^1([0,\tau_1];C^1_c(\R))$, $\zeta\ge 0$, $\zeta(\cdot,\tau_1)=0$ in $\R$, such that 
$$
\text{${\rm supp}\,\zeta\cap \big( \{x_j\}\times (0,\tau_1)\big) \neq\emptyset$\;\; for some $j=1,\dots, p$\,.}
$$ 

We only give the proof when $\zeta(x,t)=\alpha(x) \beta(t)$ with  $\alpha\in C^1_c(\R)$, $\alpha\ge0$, $\alpha(x_j)>0$  for a unique $j\in\{1,\dots, p\}$, and $\beta\in C^1([0,\tau_1])$, $\beta\ge0$, $\beta(\tau_1)=0$ (the general case can be dealt with similarly). We also assume $c_j>0$, since the proof is similar for $c_j<0$. Let us first  prove \eqref{ewf} in this case, namely
\begin{eqnarray}\label{emo}
&&\int_0^{\tau_1}\!\!\!\int_{I_j\cup I_{j+1}}
 \big[u_r\zeta_t+ H(u_r) \,\zeta_x\big]\,dxdt+ \int_{I_j\cup I_{j+1}} u_{0r}(x)\zeta(x,0)\,dx\,=\\
&&\qquad -\!\!\int_0^{\tau}\lla u_s(\cdot,t), \zeta_t(\cdot,t)\rra_{(x_{j-1},x_{j+1})}dt - \lla u_{0s},\zeta(\cdot,0)\rra_{(x_{j-1},x_{j+1})}\, \nonumber
\end{eqnarray}
for all $\zeta$ as above.  From \eqref{C_j+} we obtain 
\begin{eqnarray}\label{oz}
&&\int_0^{\tau_1}\lla u_s(\cdot,t), \zeta_t(\cdot,t)\rra_{(x_{j-1},x_{j+1})}dt+\lla u_{0s},\zeta(\cdot,0)\rra_{(x_{j-1},x_{j+1})}\,= \\
&&\quad =\alpha(x_j)\left(\int_0^{\tau_1}\!\beta'(t)C_j^+(t)\,dt \!+\!c_j \beta(0)\! \right) \!
=\!\alpha(x_j)\!\int_0^{\tau_1}\!\left[f_{x_j^+}^+(t) \!-\! f_{x_j^-}^+(t)  \right]\beta(t)\,dt. \nonumber
\end{eqnarray}

On the other hand, since  $u$ is a  solution of $(CL)$ in $Q_{j+1,\tau_1}$, by \eqref{ewf} there holds
\begin{equation*}
\iint_{Q_{j+1,\t_1}}\!\!\!
 \big\{ (u_r-k)\xi_t+ [H(u_r)-H(k)] \,\xi_x\big\}\,dxdt= - \int_{I_{j+1}}[ u_{0r}(x)-k]\,\xi(x,0)\,dx
\end{equation*}
for all $k\in\R$  and $\xi\in C^1([0,\tau_1];C^1_c(I_{j+1}))$, $\xi(\cdot,\tau_1)=0$ in $I_{j+1}$. Let  $\eta_q$ be defined by 
$$ 
\eta_q(x)=\left[2q(x-x_j)-1\right]\chi_{ \left[x_j+\frac{1}{2q},x_j+\frac{1}{q} \right]}(x)+\chi_{ \left(x_j+\frac{1}{q},x_{j+1} \right]}(x)\,,
$$
and let $\zeta\in C^1([0,\tau_1];C^1_c([x_j,x_{j+1})))$, $\zeta(\cdot,\tau_1)=0$ in $I_{j+1}$ (here $x_{j+1}=\infty$ if $j=p$). By standard arguments we can choose $\xi=\zeta\eta_q$ in the above equality. Then we get
$$
\begin{aligned}
&\iint_{Q_{j+1,\t_1}}  \!\!\!\!\big\{ (u_r-k)\zeta_t\eta_q+ [H(u_r)-H(k)] \,\zeta_x\eta_q\big\}dxdt+\\ 
&\qquad+\!\! \int_{I_{j+1}}\!\!\![ u_{0r}(x)-k]\,\zeta(x,0)\eta_q(x)dx
=\; -2q\int_0^{\tau_1}\!\!\!\int_{x_j+1/2q}^{x_j+1/q} [H(u_r)-H(k)]\,\zeta\,dxdt\,.
\end{aligned}
$$
Letting $q\to \infty$ in the above equality plainly gives (see \eqref{z1} and \eqref{defur}): 
\begin{eqnarray}\label{mkjn1}
&&\iint_{Q_{j+1,\t_1}}\!\!\!\!  \left\{ (u_r\!-\!k)\zeta_t+ [H(u_r)\!-\!H(k)] \zeta_x\right\}\,dxdt
\!+\! \int_{I_{j+1}}[ u_{0r}(x)\!-\!k]\zeta(x,0)\,dx\,=\\
&&\qquad = -\;{\rm ess}\lim_{x\to x_j^+}\int_0^{\tau_1}\![H(u_r(x,t))-H(k)] \,\zeta(x,t)\,dt=\nonumber \\
&&\qquad = -\int_0^{\tau_1}\![f_{x_j^+}^+(t)-H(k)] \,\zeta(x_j,t)\,dt\,.  \nonumber
\end{eqnarray}

Since  $u$ is an entropy solution of $(CL)$ in $Q_{j+1,\tau_1}$, arguing as before we obtain
$$
\begin{aligned}
&\iint_{Q_{j+1,\t_1}}\!\!\!\!\!\!\left\{|u_r\!-\!k|\zeta_t\!+\!\sgn(u_r\!-\!k)\left [H(u_r)\!-\!H(k)\right ]\zeta_x\right\}\!dxdt 
+ \!\!\!\int_{I_{j+1}}\!\!\!\!|u_{0r}(x)\!-\!k|\zeta(x,0)dx\! \ge \\
&\qquad\ge -\,{\rm ess}\lim_{x\to x_j^+}\int_0^{\tau_1}\sgn(u_r(x,t)-k)\left[H(u_r(x,t))-H(k)\right ]\zeta(x,t)\,dt
\end{aligned}
$$
for all  $\zeta$ as above,  $\zeta\ge0$. Choosing $\zeta(x,t)=\alpha(x) \beta(t)$ with  $\alpha\in C^1_c([x_j,x_{j+1}))$, $\beta\in C^1([0,\tau_1])$, $\a\ge0,\, \beta\ge0$ and $\beta(\tau_1)=0$, by \eqref{superlim1++} from the above inequality we obtain
\begin{eqnarray}\label{mkjn3}
&&\iint_{Q_{j+1,\t_1}}\!\!\! \left\{|u_r-k|\,\zeta_t+\sgn(u_r-k)\left [H(u_r)-H(k)\right ]\zeta_x\right\}dxdt \,+ \\
&&\quad+\! \int_{I_{j+1}}\!\!|u_{0r}(x)\!-\!k|\zeta(x,0)\,dx + {\rm ess}\!\!\lim_{x\to x_j^+}\int_0^{\tau_1}\! \!\left[H(u_r(x,t))\!-\!H(k)\right ]\zeta(x,t)\,dt\ge  \nonumber  \\
&&\quad\ge -2\,{\rm ess}\lim_{x\to x_j^+}\int_0^{\tau_1}\sgn_-(u_r(x,t)-k)\left[H(u_r(x,t))-H(k)\right ]\zeta(x,t)\,dt\,=\nonumber \\
&&\quad=-2\alpha(x_j){\rm ess}\lim_{x\to x_j^+}\int_0^{\tau_1}\sgn_-(u_r(x,t)\!-\!k)\left[H(u_r(x,t))\!-\!H(k)\right ]\beta(t)\,dt\ge 0, \nonumber
\end{eqnarray}
since $\sgn(u)=1+2\,\sgn_-(u)$. 

Replacing $Q_{j+1,\tau_1}$ by $Q_{j,\tau_1}$ we obtain, similarly to  \eqref{mkjn1}-\eqref{mkjn3},
\begin{eqnarray}\label{mkjn2}
&& \iint_{Q_{j,\t_1}}\!\!\!
 \big\{ (u_r-k)\zeta_t + [H(u_r)-H(k)] \,\zeta_x \big\}\,dxdt+ \int_{I_j}[ u_{0r}(x)-k]\,\zeta(x,0)\,dx\,=\\
&&\quad= \!{\rm ess}\!\!\lim_{x\to x_j^-}\!\!\int_0^{\tau_1}\![H(u_r(x,t))\!-\!H(k)] \zeta(x,t)\,dt
\!=\! \int_0^{\tau_1}\![f_{x_j^-}^+(t)\!-\!H(k)] \zeta(x_j,t)\,dt,  \nonumber
\end{eqnarray}
\begin{eqnarray}\label{mkjn4}
&&\ \iint_{Q_{j,\t_1}}\!\!\! \left\{|u_r-k|\,\zeta_t+\sgn(u_r-k)\left [H(u_r)-H(k)\right ]\zeta_x\right\}dxdt \,+ \\
&&\quad+\! \int_{I_j}\!\!|u_{0r}(x)\!-\!k|\zeta(x,0)\,dx - {\rm ess}\!\!\lim_{x\to x_j^-}\int_0^{\tau_1}\!\! \left[H(u_r(x,t))\!-\!H(k)\right ]\zeta(x,t)\,dt\ge 0, \nonumber 
\end{eqnarray}

Summing \eqref{mkjn1} and \eqref{mkjn2}  gives
\begin{eqnarray}\label{emo1}
&&\int_0^{\tau_1}\!\!\!\int_{I_j\cup I_{j+1}}
 \big[u_r\zeta_t+ H(u_r) \,\zeta_x\big]\,dxdt+ \int_{I_j\cup I_{j+1}} u_{0r}(x)\zeta(x,0)\,dx\,=\\
&&\qquad=-\,\alpha(x_j)\int_0^{\tau_1}\left[ f_{x_j^+}^+(t) - f_{x_j^-}^+(t) \right]\beta(t)\,dt\,.\nonumber
\end{eqnarray}
Then equality \eqref{emo} follows from \eqref{oz} and \eqref{emo1}.

Next we  prove \eqref{mkru} for all $\zeta$ as above, namely 
\begin{eqnarray}\label{mkrum}
&&\int_0^{\tau_1}\!\!\!\int
_{I_j\cup I_{j+1}}
 \left\{|u_r-k|\,\zeta_t+\sgn(u_r-k)\left [H(u_r)-H(k)\right ]\zeta_x\right\}dxdt \,+\\
&&\qquad+ \int_0^{\tau_1}\left\langle |u_s(\cdot,t)|,\zeta_t(\cdot,t)\right\rangle_{(x_{j-1},x_{j+1})}\,dt +  \left\langle |u_{0s}|, \zeta(\cdot,0)\right\rangle_{(x_{j-1},x_{j+1})} \,\ge \nonumber\\
&&\qquad \ge-\int_{I_j\cup I_{j+1}}
|u_{0r}(x)-k|\,\zeta(x,0)\,dx \nonumber
\end{eqnarray}
Since  $u$ is  an entropy solution of $(CL)$ in $Q_{j,\tau_1}$ and $Q_{j+1,\tau_1}$, from  \eqref{z1},  \eqref{z2}, \eqref{mkjn3} and \eqref{mkjn4} it follows that
\begin{eqnarray*}
&&\int_0^{\tau_1}\!\!\!\int
_{I_j\cup I_{j+1}}\!\! \left\{|u_r-k|\
\,\zeta_t+\sgn(u_r-k)\left [H(u_r)-H(k)\right ]\zeta_x\right\}dxdt \,+ \\
&&\qquad+\!  \int_{I_j\cup I_{j+1}}\!\!|u_{0r}(x)-k|\,\zeta(x,0)\,dx\, \ge\,
-\,\alpha(x_j)\int_0^{\tau_1}\left[ f_{x_j^+}^+(t) - f_{x_j^-}^+(t) \right]\beta(t)\,dt\,. \nonumber
\end{eqnarray*}
Since $|u_{0s}|\lefthalfcup (x_{j-1},x_{j+1})=u_{0s}\lefthalfcup (x_{j-1},x_{j+1})$ (recall that $c_j>0$ by assumption) and by \eqref{defus} $|u_s(\cdot,t)|\lefthalfcup (x_{j-1},x_{j+1}) =u_s(\cdot,t)\lefthalfcup (x_{j-1},x_{j+1})$ for all $t\in [0,T ]$, the above inequality together with \eqref{oz} implies \eqref{mkrum}. Therefore, 
the measure $u$ defined by \eqref{defuu} is an entropy solution of $(CL)$ in $S_{\tau_1}$. 

If $\tau_1<T$, either $u_s(\cdot,\tau_1)=0$, or $u_s(\cdot,\tau_1)\neq0$. If $u_s(\cdot,\tau_1)=0$, there holds $C_j^\pm(\tau_1)=0$ for all $j=1,\dots,p$, thus $ u_s(\cdot,t)=0$ for all $t\in[\tau_1,T]$. Then, by the standard theory of scalar conservation laws, we can continue the solution of $(CL)$ in $(\tau_1,T]$ with initial data $u_r(\cdot,\tau_1)$. On the other hand,  if $ u_s(\cdot,\tau_1)\neq0$, then $C_j^\pm(\tau_1)\neq0$ for some $j=1,\dots,p$ and, arguing as before, we can continue the solution of $(CL)$ in $(\tau_1,\tau_2]$, 
with initial data $u(\cdot,\tau_1)$,  for some $\tau_2\in (\tau_1,T]$. Iterating the procedure $q$ times with  $2\le q\le p$, we obtain that either $\tau_q=T$, or $u_s(\cdot,\tau_q)=0$. 

\smallskip

Let us now address uniqueness. 
Let $u, v\in C([0,T];\mathcal{M}(\R))$ be entropy solutions of $(CL)$, and let
$\tau:= \min\, \{t_u,t_v\}$, where
\begin{equation*}
\begin{cases}
\;t_u:= \sup\, \{t\in[0,T) \,|\,  {\rm supp}\,u_s(\cdot,t)={\rm supp}\,u_{0s} \}\, \smallskip\\
\;t_v:= \sup\, \{t\in[0,T) \,|\,  {\rm supp}\,v_s(\cdot,t)={\rm supp}\,u_{0s} \}\,.
\end{cases}
\end{equation*}
Arguing as at the end of the existence  proof, it is enough to show that $u=v $ in $\mathcal{M}(S_{\tau})$. We claim that this follows, if we prove that 
\begin{equation}\label{contra}
\text{ $u_r=v_r$ \quad a.e.~in $S_{\tau}$\,.}
\end{equation}
In fact, equalities \eqref{ewf} and \eqref{contra} imply that, for all $\zeta\in C^1([0,\tau];C^1_c(\R))$, $\zeta(\cdot,\tau)=0$ in $\R$,
$$
\int_0^{\tau}\!\! \langle u_s(\cdot,t) \! - \! v_s(\cdot,t), \zeta_t(\cdot,t)\rangle_{\R}\,dt \!
= \! \! \iint_{S_{\tau}} \! \{(u_r \! - \! v_r)\,\zeta_t+ [H(u_r)\! - \! H(u_r)] \zeta_x\}\,dxdt= 0\,.
$$
Hence  $\langle u_s(\cdot,t)-v_s(\cdot,t),\alpha\rangle_{\R}=0$ for a.e.~$t\in (0,\tau)$, for all $\alpha \in C^1_c(\R)$.
Therefore $u_s=v_s$ in $L^{\infty}(0,\tau;\mathcal{M}(\R))$, thus \eqref{contra} implies $u=v $ in $\mathcal{M}(S_{\tau})$.

It remains to prove \eqref{contra}, which is equivalent to showing that $u_r=v_r$ a.e.~in $Q_{j,\tau}$ for all $j=1,\dots,p+1$. However, this follows from the uniqueness results provided by Theorem \ref{cmp}.. Then the result follows.
\hfill$\square$

\bigskip
\noindent
{\bf Acknowledgement.} 
M.B. acknowledges the MIUR Excellence Department Project awarded to the
Department of Mathematics, University of Rome Tor Vergata, CUP E83C18000100006.


\end{document}